\documentclass[11pt]{article}
\usepackage{amsmath,amsfonts,amssymb}
\usepackage{bbm}
\usepackage{cases}
\usepackage[hidelinks]{hyperref}

\def\<{{\langle}}
\def\>{{\rangle}}
\allowdisplaybreaks[4]

\newtheorem{theorem}{Theorem}[section]
\newtheorem{lemma}[theorem]{Lemma}
\newtheorem{corollary}[theorem]{Corollary}
\newtheorem{proposition}[theorem]{Proposition}

\newtheorem{remark}[theorem]{Remark}
\newtheorem{example}[theorem]{Example}
\newtheorem{definition}[theorem]{Definition}

\begin{document}
	
	\title{\vspace*{-1.5cm}
		Large deviations of fully local monotone stochastic partial differential equations driven by gradient-dependent noise}
	
	\author{Tianyi Pan$^{1}$, Shijie Shang$^{1}$, Jianliang Zhai$^{1}$, Tusheng Zhang$^{2}$}
	\footnotetext[1]{\, School of Mathematics, University of Science and Technology of China, Hefei, China. Email:pty0512@mail.ustc.edu.cn(Tianyi Pan), sjshang@ustc.edu.cn(Shijie Shang), zhaijl@ustc.edu.cn(Jianliang Zhai).}
	\footnotetext[2]{\, School of Mathematics, University of Manchester, Oxford Road, Manchester
	M13 9PL, England, U.K. Email: tusheng.zhang@manchester.ac.uk.}
	\date{}
	\maketitle

	\begin{abstract}
	Consider stochastic partial differential equations (SPDEs) with fully local monotone coefficients in a Gelfand triple $V\subseteq H\subseteq V^*$
	\begin{numcases}{}
		du^\varepsilon_t=A(t,u^\varepsilon_t)dt+\varepsilon B(t,u^\varepsilon_t)dW_t,\ t\in (0,T],\\
		u^\varepsilon_0=x\in H,\nonumber
	\end{numcases}
	where
	$$A: [0,T] \times V\rightarrow V^*,\ \ B:[0,T]\times V\rightarrow\ L_2(U,H)$$ are measurable maps, $L_2(U,H)$ is the space of Hilbert-Schmidt operators from $U$ to $H$ and $W$ is a $U$-cylindrical Wiener process.\par
	In this paper, we establish a small noise large deviation principle(LDP) for the solutions \{$u^\varepsilon$\}$_{\varepsilon>0}$ of the above SPDEs.
	The main contribution of this paper is the much more generality of our framework than that of  the existing results. In particular, the diffusion coefficient $B(t,\cdot)$ may depend on the gradient of the solutions, which is of great interest in the field of SPDEs, but there are few existing results on the topic of LDP. The
	broader scope of the fully local monotone setting  leads us to use different strategies and techniques. A combination of the pseudomonotone technique and compactness arguement plays a crucial role in the whole paper.
	
	Our framework is very general to include many interesting models that could not be covered by existing work, including
	stochastic quasilinear SPDEs, stochastic convection diffusion equation, stochastic 2D Liquid crystal equation, stochastic $p$-Laplace equation with gradient-dependent noise, stochastic 2D Navier-Stokes equation with gradient-dependent noise etc.
\end{abstract}


\medskip

\newpage
\tableofcontents
\section{Introduction}
\ \ \ Let $H$ be a separable Hilbert space with inner product $( \cdot,\cdot )$ and norm $\|\cdot\|_H$. Let $V$ be a reflexive Banach space that is continuously, densely and compactly embedded into $H$ and let $V^*$ be the dual space of $V$. We denote the norms of $V$ and $V^*$ by $\|\cdot\|_V$ and $\|\cdot\|_{V^*}$ respectively. By the Riesz representation, the Hilbert space $H$ can be identified with its dual space $H^*$. Thus we obtain a Gelfand triple
$$V\subseteq H\subseteq V^*.$$
We denote the dual pair between $f\in V^*$ and $v\in V$ by $\langle f,v \rangle$. It is easy to see that $(u,v)=\langle u,v\rangle$, $\forall u\in H, v\in V$.\par
Consider a cylindrical Wiener process $W$ on a separable Hilbert sapce $U$ defined on some probability space $(\Omega,\mathcal{F},\mathcal{F}_t,P)$ satisfying the usual conditions.
One can interpret $W$ as an $U_1$-valued Wiener process in the larger Hilbert space $U_1$, where the embedding $U\subseteq U_1$ is Hilbert-Schimdt. Moreover, the Wiener process $W$ has continuous trajectories. \par
Let $T>0$ be fixed in this paper. Consider the following small noise perturbation of partial differential equations: for $\varepsilon>0$,
\begin{numcases}{}\label{spde}
	du^\varepsilon_t=A(t,u^\varepsilon_t)dt+\varepsilon B(t,u^\varepsilon_t)dW_t,\ t\in (0,T],\\
	u^\varepsilon_0=x\in H,\nonumber
\end{numcases}
where the maps
$$A: [0,T] \times V\rightarrow V^*,\ \ B:[0,T]\times V\rightarrow\ L_2(U,H)$$ are measurable. $L_2(U,H)$ is the space of Hilbert-Schmidt operators from $U$ to $H$, whose norm is written as $\|\cdot\|_{L_2}$ for short. \par
There exists now a rich theory on the stochastic partial differential equations(SPDEs) of the form (\ref{spde})
with monotone coefficients $A(t,\cdot)$, which means that $A(t,\cdot)$ satisfies for each $u,v\in V$,
\begin{align}\label{10120008}
	\langle A(t,u)-A(t,v), u-v\rangle \leq [f(t)+\rho(u)+\eta(v)]\|u-v\|_H^2,
\end{align}
where $\rho(\cdot)$ or $\eta(\cdot)$ are locally bounded functions on $V$ and $f:[0,T]\rightarrow\mathbb{R}_+$ is a measurable function.
The well-posedness of such SPDEs with multiplicative noise in the strictly monotone case, that is $\rho(\cdot)=\eta(\cdot)=0$, was first established by Pardoux in \cite{P} and \cite{P2} in the later 1970s. Significant contributions have been made in  subsequent works, such as \cite{KR2,G,LR,LR2,BLZ,L2,NS,NTT,KM}, regarding the well-posedness of the solution of SPDEs with monotone coefficients.  \par

Among these works, we mention that the solvability of the local monotone SPDE, where either $\eta(\cdot)=0$ or $\rho(\cdot)=0$ in (\ref{10120008}), was obtained by Liu and R\"{o}ckner in \cite{LR2}. However, the well-posedness of the fully local monotone SPDE, where both $\rho(\cdot)$ and $\eta(\cdot)$ are  nonzero, has only been established  
in a very recent paper \cite{RSZ}.
We like to mention that the framework for the well-posedness of SPDEs in \cite{RSZ} is 
considerably more general than the ones in the literature.
Many important examples are listed in \cite{RSZ}, which have precisely such a fully local monotone coefficients, rather than merely local monotone; see also
Section \ref{sec 5} in this paper.

The purpose of this paper is to establish large deviation principle(LDP) of SPDEs with fully local monotone coefficient within the very general setting presented in \cite{RSZ}, where the diffusion coefficient is also allowed to be gradient-dependent.
The large deviation principle describes the limiting behavior of the laws of solutions as the noise in the equation converges to zero, in terms of a rate function. It has wide applications in various fields, including hydrodynamics,
statistical mechanics and risk management. For a more comprehensive understanding, we refer the readers to the monograph  \cite{DZ}.

There has been a great amount of  literatures on small noise large deviations for both SDEs and SPDEs. Here we only mention some results closely related to the results in this paper.
In  \cite{RZ,L}, the authors established the LDP for SPDEs with strict monotonicity.
Obviously, the LDP for local monotone SPDEs has been considered in  \cite{CM} and  \cite{XZ,LTZ}.
The main contribution of this paper is that we cover SPDEs with the fully local monotone coefficient, which is much more general than the existing results. On the other hand, we stress that none of these previous works mentioned above for LDP allow gradient-dependent noise, that is, the diffusion coefficient $B(t,\cdot)$ depends on elements in the $V$-norm. On the other hand, many works concerning the well-posedness of  monotone SPDEs do allow gradient-dependent noise. SPDEs with gradient-dependent noise have been a crucial part of the works by Flandoli, Krylov, Rozovski\u{ı} and Brze\'{z}niak, as demonstrated in \cite{BCF,GK, MR, MR2}. Additionally, there have been some rigorous justifications of such noise in fluid equations recently, which can be found in  \cite{FGP,H2,BFM,LC,FL,FGL}.
Therefore, the LDP for SPDEs with gradient-dependent noise is very much needed to be investigated. There are only a handful of papers so far dedicated to this topic. \cite{SS} examines the LDP for stochastic Navier-Stokes equation, \cite{S} focuses on the stochastic primitive equation, \cite{MSS} investigates the LDP for the inviscid shell model and \cite{GL} explores the LDP related to the scaling limit of the stochastic Euler equation with transport noise.
Since the techniques employed in these papers rely on the specific nature of the fluid dynamics, they fail in our much general setting. This strongly motivates our investigation.
We emphasize that the conditions for the dependence on the gradient of the solutions of the diffusion coefficient we make, see (\hyperlink{H2*}{H2*})(\hyperlink{H5*}{H5*}), are weaker than the ones imposed
in the above mentioned papers.
In order to show that the solutions of SPDEs (\ref{spde}) satisfy the LDP, we will adopt the weak convergence approach introduced in \cite{BD}. More precisely,  the sufficient criteria presented in \cite{MSZ} will play an important role; see Theorem \ref{criterion}. The main ingredients of our proofs are to establish the well-posedness of solutions
to the so called skeleton equations, the convergence of the solutions of the skeleton equations  when the driving signals converge weakly, and the 
convergence of the solutions of the controlled SPDEs to  the solutions of the skeleton
equations.

\par Because of the different techniques involved, we need to distinguish two cases. In Part A, we treat the first case where the diffusion coefficient $B(t,\cdot)$ is assumed to be continuous on the space $H$, see (\hyperlink{H5}{H5}). We point out that this condition is weaker than the Lipschitz condition usually imposed in the previous works. The situation where the diffusion coefficient $B(t,\cdot)$ is gradient-dependent is handled in Part B.
The generality of the fully local monotone setting makes the identification of the limits of the solutions of the approximating equations much more challenging. 

To overcome these challenges, a combination of pseudomonotone techniques and the compactness arguments plays a crucial role, as shown in the proofs of the existence of the solution to skeleton equations and the verification of condition (b) of Theorem \ref{criterion} in both parts. Now, let us briefly describe the main technical differences between the two parts. In Part A, to carry out the weak convergence approach, we prove the compactness of the solutions of the Galerkin approximating equations and fully utilize the continuity of the diffusion coefficient and the Bochner-pseudomonotonicity of the monotone coefficient(a concept introduced in Definition 3.1 of \cite{KR}). 
On the other hand, in Part B, to address the issue arising from the lack of continuity of the diffusion coefficient and the dependence on the gradient of
the solutions, we establish the compactness of the approximating solutions and employ a modified monotone argument, which is inspired by the technique previously used in Lemma 3.4 of \cite{RSZ}. This argument is a modified version of the proof of (4.53) in \cite{LR}, which is no longer applicable in our much general case.\par 



Finally, we mention some new important examples to which our LDP results apply. These include  quasilinear stochastic partial differential equations, the convection diffusion equation, two-dimensional liquid crystal system, Cahn-Hilliard equation, two-dimensional Allen-Cahn-Navier-Stokes equation and two-dimensional Cahn-Hilliard-Navier-Stokes equation.
 These examples are covered as applications of Part A. The two-dimensional Navier-Stokes equation and $p$-Laplace equation with gradient dependent noise, will also be included as applications of Part B. These examples represent important and complex models in various scientific and engineering disciplines.

The rest of the paper is organized as follows. Section 2 provides a brief overview of the weak convergence method in the theory of LDP.  In Subsection 3.1 we introduce the hypotheses on the coefficients and state the main result in Part A. In subsection 3.2, we establish the well-posedness of the skeleton equation. Subsection 3.2 and 3.3 are devoted to the proof the main result in Part A.  Section 4 follows the same structure as Section 3, focusing on the LDP of SPDEs whose diffusion coefficients are allowed to be gradient-dependent. As applications, examples are provided in Section 5.\par
Here are some conventions. Throughout, $C$
denotes a generic positive constant whose value may vary from line to line.
Other constants will be denoted as $C_1$, $C_2$,$\cdots$. They are all positive but
their values are not important. The dependence of constants on parameters
if necessary will be indicated, e.g. $C_T$ .
	\section{Large Deviation Principle}

In this section, we recall the weak convergence method of LDP.
	\begin{definition}
		Let $\mathcal{E}$ be a Polish space with metric $\rho$ and the Borel $\sigma$-field $\mathcal{B}(\mathcal{E})$. A function I: $\mathcal{E}\rightarrow[0,+\infty]$ is called a \textbf{rate function} if the level set $\{e\in \mathcal{E}:I(e)\leq M\}$ is a compact subset of $\mathcal{E}$ for each $M<\infty$.
	\end{definition}
	
	\begin{definition}
		A family of $\mathcal{E}$-valued random variables $\{X^\varepsilon\}_{\varepsilon>0}$ is said to satisfy the LDP on $\mathcal{E}$ with rate function I if for each Borel subset B of $\mathcal{E}$,
		$$-\inf\limits_{e\in\mathring{B}}I(e)\leq \liminf\limits_{\varepsilon\rightarrow 0}\varepsilon^2 \log P(X^\varepsilon\in B)\leq\limsup\limits_{\varepsilon\rightarrow 0}\varepsilon^2 \log P(X^\varepsilon\in B)\leq -\inf\limits_{e\in\bar{B}}I(e),$$
		where $\bar{B}$ and $\mathring{B}$ denotes respectively the closure and interior of the set $B$.
	\end{definition}
	\vskip 0.25cm
	\ \ \ Next we introduce a sufficient condition for a sequence of Wiener functionals to satisfy the LDP. Let $\{W_t\}_{t\geq 0}$ be a $U$-cylindrical Wiener process defined on a filtered probability space $(\Omega,\mathcal{F},\mathcal{F}_t,P)$ as mentioned in the introduction. \par
	For any $N>0$, we define the space of $U$-valued control functions as follows:
	\begin{align*}
			S_N&\dot{=}\Big\{\phi\in L^2([0,T],U):\int_{0}^{T}\|\phi(s)\|_U^2ds \leq N\Big\},
	\end{align*}	

	We also introduce the spaces of random control processes,
	\begin{align*}
	\mathcal{A}\ \dot{=}\ \ \Big\{v:&[0,T]\times\Omega\rightarrow U, v \text{ is an $U$-valued  $\{\mathcal{F}_t\}$-predictable process such that}\\&
\int_{0}^{T}\|v(s,\omega)\|_U^2ds<\infty,\  P-a.e.\Big\},\\
		\mathcal{A}_N\ \dot{=}\ \Big\{v&\in \mathcal{A}:v(\cdot,\omega)\in S_N,\  P-a.e.\Big\}.
	\end{align*}
We remark that
 $S_N$ is a compact Polish space equipped with the weak topology of $L^2([0,T],U)$.  Let $U_1$ be a Hilbert space such that the embedding $U\subset U_1$ is Hilbert-Schmidt. Let $\mathcal{G}^\varepsilon:C([0,T],U_1)\rightarrow \mathcal{E}$ be a measurable map for each $\varepsilon>0$, and define $u^\varepsilon\dot{=}\ \mathcal{G}^\varepsilon(W)$.
	To establish the LDP of the Wiener functionals $\{u^\varepsilon\}_{\varepsilon>0}$, we will employ the following sufficient conditions established in \cite{MSZ}, which are based on a criteria of Budhiraja-Dupuis in \cite{BD}.
	
	\begin{theorem}{}\label{criterion}
		If there exists a measurable map $\mathcal{G}^0:C([0,T],U_1)\rightarrow \mathcal{E}$ such that the following two conditions hold,\\
		(a) For every $N<\infty$, for any family $\{h^\varepsilon\}_{\varepsilon>0}\subseteq \mathcal{A}_N$ and any $\delta>0$,$$\lim\limits_{\varepsilon\rightarrow 0}P\big(\rho(Y^{\varepsilon},X^{\varepsilon})>\delta\big)=0,$$
		where $X^{\varepsilon}\dot{=}\ \mathcal{G}^\varepsilon\big(W.+\frac{1}{\varepsilon}\int_{0}^{.}{h}^\varepsilon(s)ds\big)$, $Y^{\varepsilon}\dot{=}\ \mathcal{G}^0\big(\int_{0}^{.}{h}^\varepsilon(s)ds\big)$ and $\rho(\cdot,\cdot)$ stands for the metric of the space $\mathcal{E}$.\\
		(b) For every $N<\infty$ and any family $\{h^n\}_{n\in\mathbb{N}}\subseteq S_N$ that converges weakly to some element $h$ in $L^2([0,T],U)$ as $n \rightarrow \infty$ ,  we have $$\mathcal{G}^0\Big(\int_{0}^{.}{h}^n(s)ds\Big)\rightarrow\mathcal{G}^0\Big(\int_{0}^{.}{h}(s)ds\Big)\ in\ \mathcal{E}.$$
		Then the family $\{u^\varepsilon\}_{\varepsilon>0}=\{\mathcal{G}^\varepsilon(W)\}_{\varepsilon>0}$ satisfies a large deviation principle with the rate function
		\begin{align*}
			I(f)=\inf\limits_{\big\{{h}\in L^2([0,T],U):f=\mathcal{G}^0(\int_{0}^{\cdot}{h}(s)ds)\big\}}\Big\{\frac{1}{2}\int_{0}^{T}\|{h}(s)\|_U^2ds\Big\},
		\end{align*}
		with the convention $\inf\{\emptyset\}=\infty$.
	\end{theorem}
\section{Part A}
\setcounter{equation}{0}
\ \ \ In this part, we study the LDP for SPDEs with fully local monotone coefficients and $H$-continuous diffusion coefficients. In Subsection 3.1 we state the conditions on the coefficients and present the main result. In Subsection 3.2 we obtain the well-posedness of the skeleton equation. The main result is proved in  Subsection 3.3 and 3.4.\par
\subsection{Hypotheses and Main Results}
\ \ \ \ For $\varepsilon>0$, we consider the following stochastic partial differential equations:
	\begin{numcases}{}\label{spde-a}
		du^\varepsilon_t=A(t,u^\varepsilon_t)dt+\varepsilon B(t,u^\varepsilon_t)dW_t,\ t\in (0,T],\\
		u^\varepsilon_0=x\in H,\nonumber
	\end{numcases}
 We introduce the following assumptions.\par
 \vskip 0.3cm
 Let $f\in L^1([0,T],\mathbb{R}_+)$ and $\alpha\in(1,\infty)$.
\vskip -0.5cm
\begin{itemize}
	\item [\hypertarget{H1}{{\bf (H1)}}] (Hemicontinuity) For $a.e.$
	 $t\in[0,T]$, the map $\lambda\in\mathbb{R}\rightarrow\big\langle A(t,u+\lambda v),x\big\rangle \in\mathbb{R}$ is continuous, for any $u,v,x\in V$.
	\item [\hypertarget{H2}{{\bf (H2)}}] (Local Monotonicity) There exist nonnegative constants $\gamma$ and $C$ such that for $ a.e.$ $t\in[0,T]$, the following inequalities hold for any $u,v\in V$,
	\begin{align*}
		&2\big\langle A(t,u)-A(t,v), u-v\big\rangle+\|B(t,u)-B(t,v)\|^2_{L_2}\\\leq
		&\big[f(t)+\rho(u)+\eta(v)\big]\|u-v\|_H^2,\\&
		|\rho(u)|+|\eta(u)|\leq C(1+\|u\|_V^\alpha)(1+\|u\|_H^\gamma),
	\end{align*}
	where $\rho$ and $\eta$ are two measurable functions from $V$ to $\mathbb{R}$.
	\item[\hypertarget{H2'}{{\bf (H2')}}] (Generalized local monotinicity) For any $R>0$, there exists a function $K_R(\cdot)\in L^1([0,T],\mathbb{R}_+)$ such that for $a.e.\ t\in [0,T]$ and any $u,v\in V$ with $\|u\|_V\vee\|v\|_V\leq R$,
	$$\big\langle A(t,u)-A(t,v),u-v\big\rangle \leq K_R(t)\|u-v\|_H^2. $$
	\item [\hypertarget{H3}{{\bf (H3)}}] (Coercivity) There exists a constant $c>0$ such that for $a.e.$ $t\in[0,T]$, the following inequality hold for any $u\in V$,
	\begin{align*}
		2\langle A(t,u),u\rangle+\|B(t,u)\|^2_{L_2}\leq f(t)\big(1+\|u\|^2_H\big)-c\|u\|^\alpha_V.
	\end{align*}
	\item [\hypertarget{H4}{{\bf (H4)}}] (Growth) There exist nonnegative constants $\beta$ and $C$ such that for $a.e.$ $t\in[0,T]$, we have for any $u\in V$,
	$$\|A(t,u)\|^\frac{\alpha}{\alpha-1}_{V^*}\leq\big(f(t)+C\|u\|_V^\alpha\big)\big(1+\|u\|_H^\beta\big).$$
	\item [\hypertarget{H5}{{\bf (H5)}}] For $a.e.$ $t\in[0,T]$, we have
	$$\|B(t,u_n)-B(t,u)\|_{L_2}\rightarrow0$$
	for any $\{u_n\}_{n\geq1},u$ in $V$ such that $u_n$ converges to $u$ in $H$ as $n\rightarrow\infty$.
	Moreover, there exists $g\in L^1([0,T],\mathbb{R_+})$ such that for $a.e.$ $t\in[0,T]$, we have
	$$\|B(t,u)\|^2_{L_2}\leq g(t)\big(1+\|u\|_H^2\big),$$ for any $u\in V$.
\end{itemize}
\begin{remark}\label{12152249}
	According to Lemma 2.15 in \cite{RSZ}, the combination of (\hyperlink{H1}{H1}) and (\hyperlink{H2'}{H2'}) implies that for $a.e.$ $t\in[0,T]$, the operator $u\in V\rightarrow A(t,u)\in V^*$ is pseudo-monotone (see Definition 2.1 in \cite{RSZ}). Therefore, with the support of (\hyperlink{H4}{H4}) and Proposition 27.7 in \cite{Z},  the map
	$u\in V\rightarrow A(t,u)\in V^*$ is demicontinuous for $a.e.$ $t\in[0,T]$. In orther words , for $u_n\rightarrow u$ in $V$, we observe that $A(t,u_n)\rightarrow A(t,u)$ weakly in $V^*$.
\end{remark}\par
\vskip 0.3cm
The following well-posedness was obtained in \cite{RSZ}.
\begin{proposition}\label{04201820}
	Under the assumptions  (\hyperlink{H1}{H1})(\hyperlink{H2'}{H2'})(\hyperlink{H3}{H3})-(\hyperlink{H5}{H5}), there exists a probabilistic weak solution $u^\varepsilon$ to equation (\ref{spde-a}) in the sense of Definition 2.5 in \cite{RSZ}, for any initial value $x\in H$ and $\varepsilon\in(0,1]$. Moreover, if (\hyperlink{H2}{H2}) holds, then pathwise uniqueness holds for solutions of equation (\ref{spde-a}).
\end{proposition}
\vskip -0.25cm
\vskip 0.5cm
\ \ \ \ To state the main result of this section, we introduce the so called skeleton equation.\par
Let $h\in L^2([0,T],U)$. Consider the following skeleton equation
\begin{numcases}{}
	dY^h_t=A(t,Y^h_t)dt+ B(t,Y^h_t)h(t)dt,\ t\in (0,T]\label{skeleton},\\
	Y^h_0=x\in H.\nonumber
\end{numcases}
\begin{definition}\label{sol}
	We say that an element $Y^h\in C([0,T],H)\cap L^\alpha([0,T],V)$ is a solution of (\ref{skeleton}) if for all $t\in[0,T]$, the following integral equation holds in $V^*$,
	$$Y^h_t=x+\int_{0}^{t}A(s,Y^h_s)ds+\int_{0}^{t}B(s,Y^h_s)h(s)ds.$$
\end{definition}
Next result gives the existence and uniqueness of solutions to (\ref{skeleton}). The proof of this proposition can be found in Subsection 3.2.
\begin{proposition}\label{WP0}
	Under the assumptions  (\hyperlink{H1}{H1})(\hyperlink{H2'}{H2'})(\hyperlink{H3}{H3})-(\hyperlink{H5}{H5}), there exists a solution $Y^h\in C([0,T],H)\cap L^\alpha([0,T],V)$ to equation (\ref{skeleton}) for any initial value $x\in H$. Furthermore, if (\hyperlink{H2}{H2}) holds, this solution is unique in $C([0,T],H)\cap L^\alpha([0,T],V)$.
\end{proposition}

Now we can state the main result of this section.
\begin{theorem}\label{ldp}
	Assume (\hyperlink{H1}{H1})-(\hyperlink{H5}{H5}) hold. The solutions $\{u^\varepsilon\}_{0<\varepsilon\leq1}$ to equation (\ref{spde}) satisfy the large deviation principle on $C([0,T],H)$, with the rate function
	\begin{align}
		I(f)=\inf\limits_\big{\{h\in L^2([0,T],U):f=Y^h\big\}}\Big\{\frac{1}{2}\int_{0}^{T}\|{h}(s)\|_U^2ds\Big\}.
	\end{align}
\end{theorem}
\noindent {\bf Proof}. To prove Theorem \ref{ldp}, we will apply the weak convergence method outlined in Theorem \ref{criterion}.\par
Firstly, we construct the maps $\mathcal{G}^\varepsilon(\cdot)$ and $\mathcal{G}^0(\cdot)$. By Proposition \ref{04201820} and the  Yamada-Watanabe theorem, for any $\varepsilon\in(0,1]$, we obtain a measurable map $\mathcal{G}^\varepsilon(\cdot):C([0,T],U_1)\rightarrow C([0,T],H)\cap L^\alpha([0,T],V)$ such that for any $U$-cylindrical Wiener process ${W}$ on a filtered probability space,  $\mathcal{G}^\varepsilon(W)$ is the unique solution to equation (\ref{spde-a}) in $L^\alpha([0,T],V)\cap C([0,T],H)$. As a corollary of Proposition \ref{WP0}, there exists a map $\mathcal{G}^0:C([0,T],U_1)\rightarrow C([0,T],H)$ such that for any $h\in L^2([0,T],U)$, $\mathcal{G}^0(\int_{0}^{\cdot}h(s)ds)$ is the unique solution $Y^h$ given in Proposition \ref{WP0}.\par
According to Theorem \ref{criterion}, the rest of the proof will be divided into two parts.\par
Part one is to verify condition (a) in Theorem \ref{criterion}. This will be done in Subsection 3.3.
Part two is to prove that condition (b) in Theorem \ref{criterion} holds. This is done in Subsection 3.4.$\hfill\blacksquare$
\subsection{Skeleton equations}

\ \ \ \ In this subsection, we will establish the well-posedness of the skeleton equation (\ref{skeleton}). We begin by proving the uniqueness, under the assumption (\hyperlink{H2}{H2}). Let $Y^1$ and $Y^2$ be two solutions to equation (\ref{skeleton}) as defined in Definition \ref{sol}, with the same initial value $x\in H$. We can express the difference between these solutions as
\begin{align*}
	Y^1_t-Y^2_t=\int_{0}^{t}\Big(A(s,Y^1_s)-A(s,Y^2_s)\Big)ds+\int_{0}^{t}\Big(B(s,Y^1_s)h(s)-B(s,Y^2_s)h(s)\Big)ds.
\end{align*}	
Let's set $$M\ \dot{=}\sup_{t\in[0,T]}\|Y^1_t\|_H\vee\|Y^2_t\|_H.$$
By the chain rule (see Theorem 3.1 in chapter 2 of \cite{P}) and (\hyperlink{H2}{H2}), we have
\begin{align*}
	&\|Y^1_t-Y^2_t\|_H^2\\ =\ &2\int_{0}^{t}\big\langle A(s,Y^1_s)-A(s,Y^2_s),\  Y^1_s-Y^2_s\big\rangle ds+2\int_{0}^{t}\Big(B(s,Y^1_s)h(s)-B(s,Y^2_s)h(s),Y^1_s-Y^2_s\Big)ds\\ \leq & \int_{0}^{t}\Big(f(s)+\rho(Y^1_s)+\eta(Y^2_s)+\|h(s)\|_U^2\Big)\|Y^1_s-Y^2_s\|_H^2ds\\ \leq & \int_{0}^{t}\Big(f(s)+C+\|h(s)\|_U^2+CM^\gamma+CM^\gamma\|Y_s^1\|_V^\alpha+CM^\gamma\|Y_s^2\|_V^\alpha\Big)\|Y^1_s-Y^2_s\|_H^2ds.
\end{align*}	
By Gronwall's inequality, we conclude that
$\|Y^1_t-Y^2_t\|_H^2=0,\ \forall\ t\in[0,T]$. Hence, the uniqueness $Y^1\equiv Y^2$ is obtained.\par
Next we will prove the existence of a solution to equation (\ref{skeleton}), using a Galerkin approximation.\par
Let $\{e_i\}_{i=1}^\infty\subseteq V$ to be an orthonormal basis of $H$.  We define $H_n$ as the $n$-dimensional subspace of $H$ spanned by $\{e_1,\dots, e_n\}$. The orthogonal projection $P_n: V^*\rightarrow H_n$ is thus defined by
\begin{align}\label{projection}
	P_n g\ \dot{=}\ \sum_{i=1}^{n}\langle g,e_i\rangle e_i, \text{ for any $g\in V^*$.}
\end{align}
For any positive integer $n\geq 1$, we consider the following parabolic equation in the finite-dimensional space  $H_n$:
\begin{numcases}{}
	dY^n_t=P_nA(t,Y^n_t)dt+P_nB(t,Y^n_t)h(t)dt,\ t\in (0,T], \nonumber\\
	Y^n_0=P_nx\in H_n.
	\label{galerkin}
\end{numcases}
We observe that a local solution to equation (\ref{galerkin}) exists since we have indicated in Remark \ref{12152249} that the map $u\in V\rightarrow A(t,u)\in V^*$ is demicontinuous for $a.e.$ $t\in[0,T]$. Therefore, the map $u\in H_n \rightarrow P_nA(t,u)+P_nB(t,u)h(t)\in H_n$ is continuous for $a.e. \  t\in[0,T]$. Then it is well-known that the local well-posedness of equation (\ref{galerkin}) holds under the conditions (\hyperlink{H1}{H1})(\hyperlink{H2'}{H2'})(\hyperlink{H3}{H3})-(\hyperlink{H5}{H5}) (see Theorem 1 in Chapter 1 of \cite{F}).
In order to obtain the global well-posedness, we first assume that a solution $Y^n$ exists on an interval $[0,T_0)$. Then, for $t< T_0$, we can apply the chain rule (see e.g. Theorem 3.1 in chapter 2 of \cite{P}) to get
\begin{align}\label{12150137}
	&\ \|Y^n_t\|^2_H\\=&\ \|P_nx\|_H^2+2\int_{0}^{t}\big\langle A(s,Y^n_s), Y^n_s\big\rangle ds+2\int_{0}^{t}\big(B(s,Y^n_s)h(s),Y^n_s\big)ds\nonumber\\\leq&\ \|x\|_H^2+\int_{0}^{t}f(s)\big(1+\|Y^n_s\|_H^2\big)dt+\int_{0}^{t}\|h(s)\|_U^2\|Y^n_s\|_H^2ds-c\int_{0}^{t}\|Y^n_s\|_V^\alpha ds\nonumber.
\end{align}
By Gronwall's inequality, there exists a constant $C>0$ independent of $n$ and $T_0$ satisfying
\begin{align}\label{12150146}
	\sup_{t\in[0,T_0)}\|Y^n_s\|_H\leq C<+\infty.
\end{align}
By Theorem 1 in Chapter 1 of \cite{F}, we can always extend $Y^n$ to $[0,T_0+\delta]$, where $\delta>0$ only depends on the constant $C$ in $(\ref{12150146})$ for fixed $n\in\mathbb{N}$. As a result, the local solution $Y^n$ to equation (\ref{galerkin}) must be global. In particular, we can derive a uniform norm estimate for $\{Y^n\}_{n\geq1}$ in $C([0,T],H)\cap L^\alpha([0,T],V)$ based on (\ref{12150137}),  which is expressed as
\begin{align}\label{estimate0}
	\sup_{n\geq1}\Big\{\sup_{t\in[0,T]}\|Y^n_t\|_H+\int_{0}^{T}\|Y^n_s\|_V^\alpha ds\Big\}<\infty.
\end{align}
To obtain the existence of a solution to equation (\ref{skeleton}), we will prove the relative compactness of the family $\{Y^n\}_{n\geq1}$ in the following lemma.
\begin{lemma}\label{compact1}
	$\{Y^n\}_{n\geq1}$ is relatively compact in $C([0,T],V^*)\cap L^2([0,T],H)$.
\end{lemma}	
\noindent {\bf Proof}. We first prove the relative compactness in $C([0,T],V^*)$. By Theorem 3.1 in \cite{J} and (\ref{estimate0}), it suffices to show that for every $e\in H_m$, $m\in \mathbb{N}$, $\{\langle Y^n, e\rangle \}_{n\in\mathbb{N}}$ is equicontinuous as a family of real-valued functions.
Note that $$\forall e\in\mathop{\cup}_{m}H_m,\  \sup_{n\in\mathbb{N}}\|P_ne\|_V<\infty.$$
For $0\leq s\leq t\leq T$,$$ Y^n_t-Y^n_s=\int_{s}^{t}P_nA(r,Y_r^n)dr+\int_{s}^{t}P_nB(r,Y^n_r)h(r)dr.$$ Based on (\ref{estimate0}), (\hyperlink{H4}{H4}) and (\hyperlink{H5}{H5}), we assert that
\begin{align*}
	\big|\langle Y^n_t-Y^n_s, e\rangle\big|&=\int_{s}^{t}\big|\big\langle A(r,Y^n_r), P_ne \big\rangle\big| dr+\int_{s}^{t}\big|\big(B(r,Y^n_r)h(r),P_ne\big)\big|dr\\&\leq \int_{s}^{t} \|A(r,Y^n_r)\|_{V^*}\|P_ne\|_V dr+\int_{s}^{t}\|B(r,Y_r^n)\|_{L_2}\|h(r)\|_U\|P_ne\|_Hdr\\&\leq C\big(\int_{s}^{t}\|A(r,Y^n_r)\|_{V^*}^\frac{\alpha}{\alpha-1}dr\big)^\frac{\alpha-1}{\alpha}(t-s)^\frac{1}{\alpha}+C\int_{s}^{t}\|B(r,Y_r^n)\|_{L_2}\|h(r)\|_Udr\\&\leq C\big(\int_{s}^{t}\big(f(r)+C\|Y^n_r\|_V^\alpha \big)dr\big)^\frac{\alpha-1}{\alpha}(t-s)^\frac{1}{\alpha}+C\int_{s}^{t}g(r)^\frac{1}{2}\|h(r)\|_Udr\\&\leq C(t-s)^\frac{1}{\alpha}+C\int_{s}^{t}g(r)^\frac{1}{2}\|h(r)\|_Udr.
\end{align*}
We observe that the right hand side tends to $0$ uniformly in $n$ as $|t-s|\rightarrow0$. Consequently we obtain the relative compactness of $\{Y^n\}_{n\geq1}$ in $C([0,T],V^*)$. \par
To establish the relative compactness of $\{Y^n\}_{n\geq1}$ in $L^2([0,T],H)$, we apply Lemma 5.1 in \cite{RSZ}, which suggests that it suffices to prove
\begin{align*}
	\lim_{\delta\rightarrow 0}\sup_n\int_{0}^{T-\delta}\|Y^n_{t+\delta}-Y^n_t\|_{H}^2dt=0.
\end{align*}	
Indeed, integrating by parts and using (\hyperlink{H3}{H3}), (\hyperlink{H4}{H4}) and (\ref{estimate0}) we have
\begin{align*}
	&\ \|Y^n_{t+\delta}-Y^n_t\|_H^2\\=&\ 2\int_{t}^{t+\delta}\big\langle A(r,Y_r^n), Y^n_r-Y^n_t\big\rangle dr+2\int_{t}^{t+\delta}\big( B(r,Y^n_r)h(r),Y^n_r-Y^n_t\big) dr\\\leq&\ 2\int_{t}^{t+\delta}\big\langle A(r,Y_r^n), Y^n_r\big\rangle dr-2\int_{t}^{t+\delta}\big\langle A(r,Y_r^n), Y^n_t\big\rangle dr\\&+\int_{t}^{t+\delta}\|B(r,Y^n_r)\|^2_{L_2} dr+\int_{t}^{t+\delta}\|h(r)\|_U^2\|Y^n_r-Y^n_t\|_H^2dr\\\leq&\ \int_{t}^{t+\delta}f(r)\big(1+\|Y^n_r\|^2_H\big)dr+\int_{t}^{t+\delta}\|h(r)\|_U^2\|Y^n_r-Y^n_t\|_H^2dr\\&+2\|Y_t^n\|_V\int_{t}^{t+\delta}\|A(r,Y_r^n)\|_{V^*}dr\\\leq& \ \int_{t}^{t+\delta}f(r)\big(1+\|Y^n_r\|^2_H\big)dr+\int_{t}^{t+\delta}\|h(r)\|_U^2\|Y^n_r-Y^n_t\|_H^2dr\\&+2\|Y_t^n\|_V\big(\int_{0}^{T}\|A(r,Y_r^n)\|^\frac{\alpha}{\alpha-1}_{V^*}dr\big)^\frac{\alpha-1}{\alpha}\delta^\frac{1}{\alpha}\\\leq\ &C \int_{t}^{t+\delta}f(r)dr+\int_{t}^{t+\delta}\|h(r)\|_U^2\|Y^n_r-Y^n_t\|_H^2dr+C\delta^\frac{1}{\alpha}\|Y_t^n\|_V.
\end{align*}
By applying Gronwall's inequality, we obtain
\begin{align*}
	\|Y^n_{t+\delta}-Y^n_t\|_H^2\leq C_{h}\big(\int_{t}^{t+\delta}f(r)dr+\delta^\frac{1}{\alpha}\|Y^n_t\|_V\big).
\end{align*}
Integrating this inequality from $0$ to $T-\delta$ and utilizing (\ref{estimate0}) , we get
\begin{align*}
	&\int_{0}^{T-\delta}\|Y^n_{t+\delta}-Y^n_t\|_H^2dt\leq C_h \int_{0}^{T-\delta}\int_{t}^{t+\delta}f(r)dr+C_h\delta^\frac{1}{\alpha},
\end{align*}
which tends to $0 $ uniformly in $n$ as $\delta\rightarrow0$. $\hfill\blacksquare$\\\par
 From (\ref{estimate0}), (\hyperlink{H4}{H4}) and (\hyperlink{H5}{H5}), we can deduce that
\begin{align}\label{02121607}
	\sup_{n\geq1}\int_{0}^{T}\|A(t,Y^n_t)\|_{V^*}^\frac{\alpha}{\alpha-1}dt+\sup_{n\geq1}\int_{0}^{T}\|P_nB(t,Y^n_t)\|_{L_2}^2dt\leq C<\infty.
\end{align}
As a consequence of (\ref{02121607}) and Lemma \ref{compact1}, there exists a subsequence (still labeled by $\{Y^n\}_{n\geq1}$), $Y\in C([0,T],V^*)\cap L^2([0,T],H)$, $\mathcal{B}\in L^2\big([0,T],L_2(U,H)\big)$ and $\mathcal{A}\in L^\frac{\alpha}{\alpha-1}([0,T],V^*)$ such that  as $n\rightarrow\infty$,
\begin{align}\label{07281632}
	&Y^n\rightarrow Y \text{ in } C([0,T],V^*)\cap L^2([0,T],H) ,\\&Y^n\rightarrow Y\text{ weakly in $L^\alpha([0,T],V)$},\\&Y^n\rightarrow Y\text{ in the weak $*$ topology of  $L^\infty([0,T],H)$},\\&A(\cdot,Y^n_\cdot)\rightarrow \mathcal{A}\text{ weakly in $L^\frac{\alpha}{\alpha-1}([0,T],V^*)$},\\&P_nB(\cdot,Y^n_\cdot)\rightarrow \mathcal{B}\text{ weakly in $L^2\big([0,T],L_2(U,H)\big)$}.\label{3.12}
\end{align}
We define
 $$\hat{Y}_t\ \dot{=}\ x+\int_{0}^{t}\mathcal{A}(s)ds+\int_{0}^{t}\mathcal{B}(s)h(s)ds,$$
  for $t\in[0,T]$. Following the same argument as in the proof of (2.65) in \cite{RSZ}, it can be shown that $\hat{Y}=Y$, which means that for $a.e.$ $t\in[0,T]$,
 $${Y}_t\ {=}\ x+\int_{0}^{t}\mathcal{A}(s)ds+\int_{0}^{t}\mathcal{B}(s)h(s)ds.$$  In particular, this implies that $Y$ is an $H$-valued continuous function(see Theorem 3.1 in Chapter 2 of \cite{P}).  In order to show that $Y$ is the solution to equation (\ref{skeleton}), it is sufficient to prove that $\mathcal{A}(\cdot)=A(\cdot,Y_\cdot)$ and $\mathcal{B}(\cdot)=B(\cdot,Y_\cdot)$. These facts are proven in the following two lemmas.
\begin{lemma}\label{12152208}
	$\mathcal{B}(\cdot)=B(\cdot,Y_\cdot), dt\text{- almost everywhere.}$
\end{lemma}
\noindent {\bf Proof}.
Since $\|Y^n-Y\|_{L^2([0,T],H)}\rightarrow0$, there exists a subsequence (still labeled by $\{Y^n\}_{n\geq0}$) such that
\begin{align}\label{12152214}
	\lim_{n\rightarrow\infty}\|Y^n_t-Y_t\|_H=0, \text{ for $a.e.$ $t\in[0,T]$.}
\end{align}
By (\ref{estimate0}) and (\hyperlink{H5}{H5}), we conclude that
\begin{align}\label{12160146}
	\lim_{n\rightarrow\infty}\int_{0}^{T}\|P_nB(t,Y^n_t)-B(t,Y_t)\|^2_{L_2(U,H)}dt=0.
\end{align}
The uniqueness of the weak limit and (\ref{12160146}) imply that $\mathcal{B}(\cdot)=B(\cdot,Y_\cdot)$.$\hfill\blacksquare$
\\\par
Recall that the map $u\in V\rightarrow A(t,u)\in V^*$ is pseudo-monotone for $a.e.\ t\in[0,T]$, as indicated in Remark \ref{12152249}. The following lemma shows that the map
$u\in L^\alpha([0,T],V)\rightarrow A(\cdot,u_\cdot)\in L^\frac{\alpha}{\alpha-1}([0,T],V^*)$ is also pseudo-monotone in some sense.
\begin{lemma}\label{12170632} For $\{u^n\}_{n\geq1}\subseteq L^\infty([0,T],H)\cap L^\alpha([0,T],V)$ and $A(\cdot,\cdot):[0,T]\times V\rightarrow V^*$ as introduced in (\hyperlink{H1}{H1})(\hyperlink{H2'}{H2'})(\hyperlink{H3}{H3})(\hyperlink{H4}{H4}),
	suppose there exist $\mathcal{A}\in L^\frac{\alpha}{\alpha-1}([0,T],V^*)$ and $u\in L^\infty([0,T],H)\cap L^\alpha([0,T],V)$ such that as $n\rightarrow\infty$,
	\begin{align}
		&u^n\rightarrow u\text{ weakly in $ L^\alpha([0,T],V)$}\nonumber,\\	&u^n\rightarrow u\text{ in the weak $\ast$  topology of $ L^\infty([0,T],H)$}\nonumber,\\&u^n\rightarrow u\text{ in $ L^2([0,T],H)$}\nonumber,\\&A(\cdot,u_\cdot^n)\rightarrow \mathcal{A}\text{ weakly in $ L^\frac{\alpha}{\alpha-1}([0,T],V^*)$}\label{12152321},\\&\liminf_{n\rightarrow\infty}\int_{0}^{T}\big\langle A(t,u^n_t), u^n_t\big\rangle dt\geq\int_{0}^{T}\big\langle \mathcal{A}(t),u_t\big\rangle dt,\label{12152329}
	\end{align}
	then $\mathcal{A}(\cdot)=A(\cdot,u_\cdot),\ dt-a.e$.
\end{lemma}
We omit the proof here and refer the readers to Lemma 2.16 in \cite{RSZ} or Proposition 3.7(i) in \cite{KR} for detailed explanation.$\hfill\blacksquare$\par
\vskip 0.5cm
	Therefore, to demonstrate that $\mathcal{A}(\cdot)=\ A(\cdot,Y_\cdot)$, it suffices to show that ($\ref{12152329}$) holds for $\{u^n\}_{n\geq1}$ and $u$ replaced by $\{Y^n\}_{n\geq1}$ and $Y$. In fact, using the chain rule, we have
\begin{align}
	&\|Y^n_T\|_H^2=\|P_nx\|_H^2+2\int_{0}^{T}\big\langle A(t,Y^n_t),Y^n_t\big\rangle dt+2\int_{0}^{T}\big(B(t,Y^n_t)h(t),Y^n_t\big)dt\label{12160135},\\&\|Y_T\|_H^2=\|x\|_H^2+2\int_{0}^{T}\big\langle \mathcal{A}(t),Y_t\big\rangle dt+2\int_{0}^{T}\big(B(t,Y_t)h(t),Y_t\big)dt\label{12160136}.
\end{align}
Since $\|Y^n-Y\|_{C([0,T],V^*)}\rightarrow0$, we can apply the lower semi-continuity of $\|\cdot\|_H$ in $V^*$ and Fatou's lemma to obtain
\begin{align}
	\|Y_T\|_H^2\leq\liminf_{n\rightarrow\infty}\|Y^n_T\|_H^2.\label{12160138}
\end{align}
If we show
\begin{align}\label{12160142}
	\int_{0}^{T}\big(B(t,Y^n_t)h(t),Y^n_t\big)dt\rightarrow\int_{0}^{T}\big(B(t,Y_t)h(t),Y_t\big)dt,
\end{align}
combining with (\ref{12160135})-(\ref{12160142}), we can obtain ($\ref{12152329}$). We note that (\ref{12160142}) can be obtained easily by (\hyperlink{H5}{H5}), $(\ref{12152214}),(\ref{12160146}),(\ref{estimate0})$ and the dominated convergence theorem. \par
We thus have finished proving Proposition \ref{WP0}.$\hfill\blacksquare$
\begin{remark}\label{08112057}
	By the argument as in the proof of (\ref{estimate0}) in this subsection, we can infer that for any $N>0$,
	\begin{align}\label{04212055}
		\sup_{h\in S_N}\Big\{\sup_{t\in[0,T]}\|Y^h_t\|^2_H+\int_{0}^{T}\|Y^h_t\|_V^\alpha dt\Big\}<\infty,
	\end{align}
where $Y^h$ is the unique solution to equation (\ref{skeleton}).
\end{remark}
\subsection{The proof of the main result: part one}
\ \ \ \ In this subsection, we will carry out the proof of part one of Theorem \ref{ldp}, namely we will verify condition (a) of Theorem \ref{criterion} for the maps $\mathcal{G}^\varepsilon$ and $\mathcal{G}^0$ as provided at the beginning of the  proof of Theorem \ref{ldp} in Subsection 3.1.
By the Girsanov's transformation, it is clear that for any ${h}^\varepsilon\in \mathcal{A}_N$ and $\varepsilon\in(0,1]$,  $X^\varepsilon=\mathcal{G}^\varepsilon\big(W_\cdot+\frac{1}{\varepsilon}\int_{0}^{\cdot}{h}^\varepsilon(s)ds\big)$ is the unique solution of the following controlled stochastic partial differential equation,
\begin{numcases}{}
	dX^\varepsilon_t=A(t,X^\varepsilon_t)dt+\varepsilon B(t,X^\varepsilon_t)dW_t+B(t,X^\varepsilon_t)h^\varepsilon(t)dt,\ t\in (0,T]\label{spde2},\\
	X^\varepsilon_0=x\in H\nonumber.
\end{numcases}
 Firstly, we will provide some moment estimates that will play a key role in the verification of condition $(a)$.\par
\begin{proposition}\label{estimate1}
	For any $q\geq 2$ , there exists a positive constant $C_{N,q}$ such that,
	\begin{align*}
		&\sup_{\varepsilon\leq1}E\Big[\sup_{t\in[0,T]}\|X^\varepsilon_t\|_H^q+\Big(\int_{0}^{T}\|X^\varepsilon_s\|^\alpha_Vds\Big)^\frac{q}{2}\Big]\leq C_{N,q}.
	\end{align*}
\end{proposition}
\noindent {\bf Proof}.
Fix $q\in[2,\infty)$. Applying $\rm It\hat{o}$'s formula (see e.g Theorem 3.1 in Section 2 of \cite{P}), (\hyperlink{H3}{H3}) , (\hyperlink{H5}{H5}) and Young's inequality, we have for $0<\varepsilon\leq 1$,
\begin{align}\label{12162022}
	&\|X^{\varepsilon}_t\|_H^q\\=&\ \|x\|_H^q+q\int_{0}^{t}\|X^\varepsilon_s\|_H^{q-2}\big\langle A(s,X^\varepsilon_s), X^{\varepsilon}_s\big\rangle ds+q\varepsilon\int_{0}^{t}\|X^\varepsilon_s\|_H^{q-2}\big(B(s,X^\varepsilon_s)dW_s,X^\varepsilon_s\big)\nonumber\\&\ +\frac{q\varepsilon^2}{2}\int_{0}^t\|X_s^\varepsilon\|_H^{q-2}\|B(s,X^\varepsilon_s)\|_{L_2}^2ds+q\int_{0}^{t}\|X^\varepsilon_s\|_H^{q-2}\big(B(s,X^\varepsilon_s)h^\varepsilon(s),X^\varepsilon_s\big)ds \nonumber \\&\ + \frac{q(q-2)\varepsilon^2}{2}\int_{0}^{t}\|X^\varepsilon_s\|_H^{q-4}\|B(s,X^\varepsilon_s)^*X^\varepsilon_s\|_U^2ds\nonumber\\\leq&\ \|x\|_H^q+C_q\int_{0}^t\big(f(s)+g(s)+\|h^\varepsilon(s)\|_U^2\big)ds-\frac{qc}{2}\int_{0}^{t}\|X^\varepsilon_s\|_H^{q-2}\|X^\varepsilon_s\|_V^\alpha ds\nonumber\\&\ +\frac{q}{2}\int_{0}^tf(s)\|X_s^\varepsilon\|_H^qds+qC\int_{0}^{t}\|X^\varepsilon_s\|_H^{q}\big(g(s)+\|h^\varepsilon(s)\|_U^2\big)ds\nonumber\\&\ + q\varepsilon\int_{0}^{t}\|X^\varepsilon_s\|_H^{q-2}\big(B(s,X^\varepsilon_s)dW_s,X^\varepsilon_s\big).
\end{align}
By Gronwall's inequality, we have
\begin{align}\label{07250212}
	&\nonumber	\sup_{s\in[0,t]}\|X^{\varepsilon}_s\|_H^q+\frac{qc}{2} \int_{0}^{t}\|X^\varepsilon_s\|_H^{q-2}\|X^\varepsilon_s\|_V^\alpha ds\\\leq&\ C_{1}\Big(C_N+q\sup_{s\leq t}\big|\int_{0}^{s}\|X^\varepsilon_r\|_H^{q-2}\big(B(r,X^\varepsilon_r)dW_r,X^\varepsilon_r\big)\big|\Big).
\end{align}
Set
\begin{align}\label{12170252}
	\tau^{\varepsilon,M}\dot{=}\inf\Big\{t\geq0:\|X^\varepsilon_t\|_H\geq M\ \text{or}\ \int_{0}^{t}\|X^\varepsilon_s\|_V^\alpha ds\geq M\Big\}\wedge T.
\end{align}
It is evident that  $\tau^{\varepsilon,M}\uparrow T$ as $M\rightarrow\infty$, since $X^\varepsilon$ is the global solution to equation (\ref{spde2}).\par
For the last term in (\ref{07250212}), by the BDG inequality and Young's inequality, we obtain
\begin{align}\label{12162058}
	&E\Big[\sup_{s\in[0,t\wedge\tau^{\varepsilon,M}]}\big|\int_{0}^{s}\|X^\varepsilon_r\|_H^{q-2}\big(B(r,X^\varepsilon_r)dW_r,X^\varepsilon_r\big)\big|\Big]\\\leq&\  CE\Big[\Big(\int_{0}^{t\wedge\tau^{\varepsilon,M}}\|X^\varepsilon_s\|_H^{2q-2}\|B(s,X^\varepsilon_s)\|_{L_2}^2ds\Big)^\frac{1}{2}\Big]\nonumber\\\leq&\ C E\Big[\Big(\sup_{s\in[0,t\wedge\tau^{\varepsilon,M}]}\|X^\varepsilon_s\|^q_H\times\int_{0}^{t\wedge\tau^{\varepsilon,M}}\|X^\varepsilon_s\|_H^{q-2}\|B(s,X^\varepsilon_s)\|_{L_2}^2ds\Big)^\frac{1}{2}\Big]\nonumber\\\leq&\ \varepsilon_0E\Big[\sup_{s\in[0,t\wedge\tau^{\varepsilon,M}]}\|X^\varepsilon_s\|^q_H\Big]+C_{\varepsilon_0}E\Big[\int_{0}^{t\wedge\tau^{\varepsilon,M}}\|X^\varepsilon_s\|_H^{q-2}\|B(s,X^\varepsilon_s)\|_{L_2}^2ds\Big].\nonumber
\end{align}
Letting $\varepsilon_0\ \dot{=}\ \frac{1}{2qC_{1}}$. By (\ref{07250212}),  (\ref{12162058}),  (\ref{12162022}) and (\hyperlink{H5}{H5}), we have
\begin{align*}
	&E\Big[\sup_{s\in[0,t\wedge\tau^{\varepsilon,M}]}\|X^{\varepsilon}_s\|_H^q\Big]+{qc}E\Big[ \int_{0}^{t\wedge\tau^{\varepsilon,M}}\|X^\varepsilon_s\|_H^{q-2}\|X^\varepsilon_s\|_V^\alpha ds\Big]\\\leq&\  C_{N,q}+C_{N,q}\int_{0}^{t}\big(1+g(s)\big)E\big[\|X_{s\wedge\tau^{\varepsilon,M}}^{\varepsilon}\|_H^q\big]ds.
\end{align*}
Again by Gronwall's inequality, we derive
\begin{align*}
	E\Big[\sup_{s\in[0,t\wedge\tau^{\varepsilon,M}]}\|X^{\varepsilon}_s\|_H^q\Big]+{qc}E\Big[ \int_{0}^{t\wedge\tau^{\varepsilon,M}}\|X^\varepsilon_s\|_H^{q-2}\|X^\varepsilon_s\|_V^\alpha ds\Big]\leq C_{N,q}.
\end{align*}
Letting $M\rightarrow\infty$ and applying Fatou's lemma, we obtain
\begin{align}\label{12162244}
	E\Big[\sup_{s\in[0,t]}\|X^{\varepsilon}_s\|_H^q\Big]+{qc}E\Big[ \int_{0}^{t}\|X^\varepsilon_s\|_H^{q-2}\|X^\varepsilon_s\|_V^\alpha ds\Big]\leq C_{N,q}.
\end{align}
Using (\hyperlink{H3}{H3}) and (\ref{12162022}) with $q=2$, we have
\begin{align}\label{12162203}
	\nonumber&\ \|X^{\varepsilon}_t\|_H^2+c\int_{0}^{t}\|X^\varepsilon_s\|_V^\alpha ds\\\leq& \ \|x\|_H^2+2\varepsilon\int_{0}^{t}\big(B(s,X^\varepsilon_s)dW_s,X^\varepsilon_s\big) \\&\ +\int_{0}^t\big(f(s)+g(s)+\|h^\varepsilon(s)\|_U^2\big)\big(1+\|X^\varepsilon_s\|_H^2\big)ds\nonumber.
\end{align}
Hence, by Gronwall's inequality, we can derive that
\begin{align*}
	&\ \|X^{\varepsilon}_t\|_H^2+c\int_{0}^{t}\|X^\varepsilon_s\|_V^\alpha ds\leq C_{N}\Big(C_{N}+2\sup_{s\in[0,t]}\big|\int_{0}^{s}\big(B(r,X_r^\varepsilon) dW_r, X^\varepsilon_r\big)\big|\Big).
\end{align*}
For $q\geq2$, utilizing the BDG inequality, (\hyperlink{H5}{H5}) and (\ref{12162244}), we obtain
\begin{align*}\label{12162214}
	&\ E\Big[\big(\int_{0}^{t}\|X^\varepsilon_s\|_V^\alpha ds\big)^\frac{q}{2}\Big]\\\ \leq& \ C_{N,q}\Big(C_{N,q}+E\Big[\sup_{r\in[0,t]}\big|\int_{0}^{r}\big(B(s,X_s^\varepsilon) dW_s, X^\varepsilon_s\big)\big|^\frac{q}{2}\Big]\Big)\nonumber\\\leq& \ C_{N,q}\Big(C_{N,q}+E\Big[\big(\int_{0}^{t}\|B(s,X_s^\varepsilon)\|_{L_2}^2\|X_s^\varepsilon\|_H^2ds\big)^\frac{q}{4}\Big]\Big)\nonumber\\\leq& \ C_{N,q}\Big(C_{N,q}+E\Big[\big(\int_{0}^{t}g(s)(1+\|X_s^\varepsilon\|_H^4)ds\big)^\frac{q}{4}\Big]\Big)\nonumber\\\leq&\ C_{N,q}.\\& \text{Combining with (\ref{12162244}) completes the proof of Theorem \ref{estimate1}.}\qquad\qquad \blacksquare
\end{align*}
\par
By applying Chebyshev's inequality, we get the following corollary
\begin{corollary}\label{12170329} As $M\rightarrow\infty$, we have
	$$\sup_{\varepsilon\leq1}P\Big(\sup_{t\in[0,T]}\|X^\varepsilon_t\|_H>M\text{ or }\int_{0}^{T}\|X^\varepsilon_s\|_V^\alpha ds>M\Big)\rightarrow0.$$
\end{corollary}\par
\vskip 0.3cm
\textbf{\textit{Proof of Theorem 3.5: part one}}\par
\vskip 0.3cm
For simplicity, we denote $Y^\varepsilon\dot{=}\ \mathcal{G}^0\big(\int_{0}^{\cdot}h^\varepsilon(s)ds\big)$.  In particular, $Y^\varepsilon$ satisfies (\ref{skeleton}) according to  the definition of $\mathcal{G}^0$. So we can express the difference between $X^\varepsilon$ and $Y^\varepsilon$ as
\begin{align*}
	X^\varepsilon_t-Y^\varepsilon_t=&\int_{0}^{t}\big(A(s,X^\varepsilon_s)-A(s,Y^\varepsilon_s)\big)ds+\int_{0}^{t}\big(B(s,X_s^\varepsilon)h^\varepsilon(s)-B(s,Y^\varepsilon_s)h^\varepsilon(s)\big)ds\\&+\varepsilon\int_{0}^{t}B(s,X_s^\varepsilon)dW_s.
\end{align*}
By $\rm It\hat{o}$'s formula, we have
\begin{align}
	&\ \|X^\varepsilon_t-Y^\varepsilon_t\|^2_H\nonumber\\=&\ 2\int_{0}^{t}\big\langle A(s,X^\varepsilon_s)-A(s,Y^\varepsilon_s), X^\varepsilon_s-Y^\varepsilon_s \big\rangle ds+2\int_{0}^{t}\big(B(s,X^\varepsilon_s)h^\varepsilon(s)-B(s,Y^\varepsilon_s)h^\varepsilon(s), X^\varepsilon_s-Y^\varepsilon_s\big)ds\nonumber\\&+2\varepsilon\int_{0}^{t}\big(B(s,X^\varepsilon_s)dW_s, X^\varepsilon_s-Y^\varepsilon_s\big)+\varepsilon^2\int_{0}^{t}\|B(s,X^\varepsilon_s)\|_{L_2}^2ds\nonumber\\\leq&\ 2\int_{0}^{t}\big\langle A(s,X^\varepsilon_s)-A(s,Y^\varepsilon_s),X^\varepsilon_s-Y^\varepsilon_s\big\rangle ds+\int_{0}^{t}\|B(s,X^\varepsilon_s)-B(s,Y^\varepsilon_s)\|_{L_2}^2ds\nonumber\\&+\int_{0}^{t}\|h^\varepsilon(s)\|_U^2\|X^\varepsilon_s-Y^\varepsilon_s\|_H^2ds+\varepsilon M^\varepsilon_t,
\end{align}
where
$$M^\varepsilon_t\ \dot{=}\ 2\sup_{r\in[0,t]}\Big|\int_{0}^{r}\big(B(s,X^\varepsilon_s)dW_s, X^\varepsilon_s-Y^\varepsilon_s\big)\Big|+\varepsilon\int_{0}^{t}\|B(s,X^\varepsilon_s)\|_{L_2}^2ds.$$
Using (\hyperlink{H2}{H2}), it follows that
\begin{align}\label{12170309}
	&\sup_{s\in[0,t]}\|X^\varepsilon_s-Y^\varepsilon_s\|_H^2\\\leq&\ \int_{0}^{t}\big(\|h^\varepsilon(s)\|_U^2+f(s)+\rho(X^\varepsilon_s)+\eta(Y^\varepsilon_s)\big)\|X^\varepsilon_s-Y^\varepsilon_s\|_H^2ds\nonumber+\varepsilon M^\varepsilon_t.
	\nonumber\end{align}
Recalling Proposition \ref{estimate1}, the definition of $\tau^{\varepsilon,M}$ in $(\ref{12170252})$ and replacing $t$ by $t\wedge\tau^{\varepsilon,M}$ in (\ref{12170309}), we have, by Gronwall's inequality and Remark \ref{08112057},
\begin{align*}
	&\sup_{s\in[0,T\wedge\tau^{\varepsilon,M}]}\|X^\varepsilon_s-Y^\varepsilon_s\|_H^2\leq\ \varepsilon C_{M,N}M^\varepsilon_{T}.
\end{align*}
From the BDG inequality, (\hyperlink{H5}{H5}), Remark \ref{08112057} and Proposition \ref{estimate1}, we obtain
\begin{align}\label{12170401}
	\sup_{\varepsilon\leq1}E[M_T^\varepsilon]<\infty.
\end{align}
Then for  $\forall\ \varepsilon>0$,
\begin{align*}
	&\ P\Big(\sup_{t\in[0,T]}\|X^\varepsilon_t-Y^\varepsilon_t\|_H>\delta\Big)\\\leq&\  P\Big(\sup_{t\in[0,T]}\|X^\varepsilon_t-Y^\varepsilon_t\|_H>\delta,\tau^{\varepsilon,M}\geq T\Big)+P\Big(\tau^{\varepsilon,M}< T\Big)\\\leq&\ \frac{\varepsilon E\big[M_T^\varepsilon\big]}{\delta^2}C_{M,N}+\sup_{\epsilon\leq1}P(\tau^{\epsilon,M}< T).
\end{align*}
We first let $\varepsilon$ tends to $0$. By estimate (\ref{12170401}), we obtain
\begin{align*}
	&\ \lim_{\varepsilon\rightarrow0} P\Big(\sup_{t\in[0,T]}\|X^\varepsilon_t-Y^\varepsilon_t\|_H>\delta\Big)\leq\sup_{\epsilon\leq1}P(\tau^{\epsilon,M}< T).
\end{align*}
Then, letting $M\rightarrow\infty$, Corollary \ref{12170329} gives
\begin{align*}
	&\ \lim_{\varepsilon\rightarrow0} P\Big(\sup_{t\in[0,T]}\|X^\varepsilon_t-Y^\varepsilon_t\|_H>\delta\Big)=0.
\end{align*}
This completes the proof of part one. $\hfill\blacksquare$
\subsection{The proof of the main result: part two}
\ \ \ \
This subsection is devoted to the verification of (b) in Theorem \ref{criterion}. Let's denote $Y^n\ \dot{=}\ \mathcal{G}^0\big(\int_{0}^{\cdot}{h}^n(s)ds\big)$ and $Y^h\ \dot{=}\ \mathcal{G}^0\big(\int_{0}^{\cdot}{h}(s)ds\big)$, where $\{h^n,h\}_{n\geq1}\subseteq S_N$. We will show that if $ h^n\rightarrow h$ weakly in $L^2([0,T],U)$ as $n\rightarrow\infty$, then
\begin{align}\label{07272230}
	\sup_{s\in[0,T]}\|Y^n_s-Y^h_s\|_H\rightarrow0.
\end{align}
By (\ref{04212055}), we have
\begin{align}\label{estimate2}
	\sup_{n\geq1}\Big\{\sup_{t\in[0,T]}\|Y^n_t\|_H+\int_{0}^{T}\|Y^n_s\|_V^\alpha ds\Big\}<\infty.
\end{align}
Using a similar approach as in the proof of Lemma \ref{compact1}, we can conclude that $\{Y^n\}_{n\geq1}$ is relatively compact in $C([0,T],V^*)\cap L^2([0,T],H)$.
Therefore, there exist $Y\in C([0,T],V^*)\cap L^2([0,T],H)$ and
a subsequence (also labeled as $\{Y^n\}_{n\geq1}$) such that as $n\rightarrow\infty$,
\begin{align}\label{07271533}
	&Y^n\rightarrow Y\text{ in $C([0,T],V^*)\cap L^2([0,T],H)$}.
\end{align}
(\ref{estimate2}) implies that
$$Y\in L^\infty([0,T],H)\cap L^\alpha([0,T],V).$$

Notice that for any $t\in[0,T]$, the difference between $Y^n$ and $Y^h$ is expressed by
\begin{align*}
	Y^n_t-Y^h_t=\int_{0}^{t}\big(A(s,Y^n_s)-A(s,Y^h_s)\big)ds+\int_{0}^{t}\big(B(s,Y^n_s)h^{n}(s)-B(s,Y^h_s)h(s)\big)ds.
\end{align*}
By the chain rule and (\hyperlink{H2}{H2}),
\begin{align*}
	&\ \|Y^h_t-Y^n_t\|_H^2\\=&\ 2\int_{0}^{t}\big\langle A(s,Y^h_s)-A(s,Y^n_s), Y^h_s-Y^n_s\big\rangle ds+2\int_{0}^{t}\big(B(s,Y^h_s)h(s)-B(s,Y^n_s)h^{n}(s),Y^h_s-Y^n_s\big)ds\\=&\ 2\int_{0}^{t}\big\langle A(s,Y^h_s)-A(s,Y^n_s), Y^h_s-Y^n_s\big\rangle ds+2\int_{0}^{t}\big(B(s,Y^h_s)h^n(s)-B(s,Y^n_s)h^{n}(s),Y^h_s-Y^n_s\big)ds\\&\ +2\int_{0}^{t}\big(B(s,Y^h_s)h(s)-B(s,Y^h_s)h^{n}(s),Y^h_s-Y^n_s\big)ds\\\leq&\ 2\int_{0}^{t}\big\langle A(s,Y^h_s)-A(s,Y^n_s), Y^h_s-Y^n_s\big\rangle ds+\int_{0}^{t}\|B(s,Y^h_s)-B(s,Y^n_s)\|_{L_2}^2ds\\& + \int_{0}^{t}\|h^{n}(s)\|_U^2\|Y^h_s-Y^n_s\|_H^2ds\ +2\int_{0}^{t}\big(B(s,Y^h_s)h(s)-B(s,Y^h_s)h^{n}(s),Y^h_s-Y^n_s\big)ds\\\leq&\ \int_{0}^{t}\big(f(s)+\|h^n(s)\|_U^2+\rho(Y^h_s)+\eta(Y^n_s)\big)\|Y^n_s-Y^h_s\|_H^2ds+\big|Q^n_t\big|,
\end{align*}
where $$Q^n_t\ \dot{=}\ 2\int_{0}^{t}\big(B(s,Y^h_s)h(s)-B(s,Y^h_s)h^n(s),Y^h_s-Y^n_s\big)ds.$$
From ($\ref{estimate2}$), (\hyperlink{H2}{H2}) and Gronwall's inequality, we obtain
\begin{align*}
	\sup_{t\in[0,T]}\|Y^n_t-Y^h_t\|_H^2\leq C\sup_{t\in[0,T]}\big|Q^n_t\big|.
\end{align*}
By ($\ref{estimate2}$), it is easy to see that $\{Q^n_t, t\in [0, T], n\geq 1\}$ is a family of equi-uniform continuous functions. To complete the proof of part two, it suffices to show that $Q^n_t\rightarrow0$ for every $t\in [0, T]$ as $n\rightarrow\infty$. Denote
\begin{align*}
	\Gamma^n_M\ \dot{=}\ \Big\{t\in[0,T]:\|h^n(t)\|_U>M\Big\}.
\end{align*}
Since $\{h^{n}\}_{n\geq1}\subseteq S_N$, the Lebesgue measure $|\Gamma^n_M|$ of $\Gamma^n_M$ can be estimated as follows
\begin{align}\label{12170728}
	\big|\Gamma^n_M\big|\leq\frac{N}{M^2}.
\end{align}
Observe that
\begin{align*}
	Q^n_t&=2\int_{0}^{t}\big(B(s,Y^{h}_s)h(s)-B(s,Y^h_s)h^n(s),Y^h_s-Y^{n}_s\big)ds\\&=	2\int_{0}^{t}\big(B(s,Y^{h}_s)h(s)-B(s,Y^h_s)h^n(s),Y^h_s-Y_s\big)ds\\
&\ \ \  +2\int_{0}^{t}\big(B(s,Y^h_s)h^n(s),Y^n_s-Y_s\big)ds+2\int_{0}^{t}\big(B(s,Y^h_s)h(s),Y_s-Y^n_s\big)ds\\&\ \dot{=}\ \text{I$^n$+II$^n$+III$^n$}.
\end{align*}
Since
$ h^n\rightarrow h$ weakly in $L^2([0,T],U)$ as $n\rightarrow\infty$
and $B(\cdot,Y^h_\cdot)^*(Y^h_\cdot-Y_\cdot)\in L^2([0,T],U)$, then we have  $$\lim_{n\rightarrow\infty}\text{I}^n= 0.$$
 For the convergence of II$^n$, we note that
\begin{align}
	&\Big|\int_{0}^{t}\big(B(s,Y^h_s)h^{n}(s),Y^{n}_s-Y_s\big)ds\Big|\nonumber\\\leq&\ CN\Big(\int_{\Gamma^n_M}g(s)ds\Big)^\frac{1}{2}+CM\int_{0}^{T}g(s)^\frac{1}{2}\|Y_s^{n}-Y_s\|_Hds.\label{08112234}
\end{align}
First letting $n\rightarrow\infty$ and then $M\rightarrow\infty$, we obtain $\text{II}^n\rightarrow0$ by (\ref{12170728}). Similarly, we also have
$\text{III}^n\rightarrow0$. Thus we have proved that $Q^n_t\rightarrow0$. Due to the uniqueness of the solutions of the skeleton equation, the limit holds for the whole sequence and thus the proof is completed.
$\hfill\blacksquare$

\section{Part B}
\ \  \ \ In this part,  we establish the LDP for  equation (\ref{spde}) with the diffusion coefficient depending on elements in $V$-norm. This means that the diffusion coefficient could depend on the gradient of the solution when it comes to the SPDEs.  We point out that Part A and Part B do not cover each other. In Subsection 4.1, we introduce the conditions on the coefficients and state the main result of this section. In Subsection 4.2, we establish the well-posedness of the skeleton equation. Finally, in Subsections 4.3 and 4.4, we prove the main result of this section.\par
\subsection{Hypotheses and Main Results}
\setcounter{equation}{0}
\ \ \ \
For $\varepsilon>0$, we again consider the following stochastic partial differential equations:
	\begin{numcases}{}\label{spde-b}
		du^\varepsilon_t=A(t,u^\varepsilon_t)dt+\varepsilon B(t,u^\varepsilon_t)dW_t,\ t\in (0,T],\\
		u^\varepsilon_0=x\in H,\nonumber
	\end{numcases}
 We introduce the following assumptions.\par
 \vskip 0.3cm
Let $f\in L^1([0,T],\mathbb{R}_+)$ and $\alpha\in(1,\infty)$, $\beta\in[0,\infty)$. \\
\begin{itemize}
	\item [\hypertarget{H1*}{{\bf (H1*)}}] (Hemicontinuity) For $a.e.$ $t\in[0,T]$, the map $\lambda\in\mathbb{R}\rightarrow\big\langle A(t,u+\lambda v),x\big\rangle \in\mathbb{R}$ is continuous, for any $u,v,x\in V$.
	\item [\hypertarget{H2*}{{\bf (H2*)}}] (Local Monotonicity) There exist  nonnegative constants $\theta\in[0,\alpha)$ , $\delta>0$, $\gamma,\kappa$ and $C$ such that for $a.e.$ $t\in[0,T]$, and any $u,v\in V$,
	\begin{align*}
		&2\big\langle A(t,u)-A(t,v), u-v\big\rangle+\delta^2\|B(t,u)-B(t,v)\|^2_{L_2}\\\leq
		&\ [f(t)+\rho(u)+\eta(v)]\|u-v\|_H^2,
	\end{align*}
	where $\rho$ and $\eta$ are two measurable functions from $V$ to $\mathbb{R}$ such that \begin{align}&
		|\rho(u)|\leq C(1+\|u\|_H^\kappa)+C\|u\|_V^\theta(1+\|u\|_H^\gamma),\\&
		|\eta(u)|\leq C(1+\|u\|_H^{2+\beta})+C\|u\|_V^\alpha(1+\|u\|_H^\beta).
	\end{align}
	\item [\hypertarget{H3*}{{\bf (H3*)}}] (Coercivity) There exists a constant $c>0$  such that for $a.e.$ $t\in[0,T]$, the following inequality hold for any $u\in V$ with some $p>1$,
	\begin{align*}
		2\big\langle A(t,u),u\big\rangle+(p-1)\|B(t,u)\|^2_{L_2}\leq f(t)(1+\|u\|^2_H)-c\|u\|^\alpha_V.
	\end{align*}
	\item [\hypertarget{H4*}{{\bf (H4*)}}] (Growth) There exist nonnegative constants $\beta\geq 0$ and $C>0$ such that for $a.e.$ $t\in[0,T]$, we have for any $u\in V$,
	$$\|A(t,u)\|^\frac{\alpha}{\alpha-1}_{V^*}\leq f(t)(1+\|u\|_H^{2+\beta})+C\|u\|_V^\alpha(1+\|u\|_H^\beta).$$
	\item [\hypertarget{H5*}{{\bf (H5*)}}] There exists $g\in L^1([0,T],\mathbb{R_+})$ and constant $L_B>0$ such that for $a.e.$ $t\in[0,T]$, we have
	$$\|B(t,u)\|^2_{L_2}\leq g(t)(1+\|u\|_H^2)+L_B(1+\|u\|_V^\alpha),$$
	for any $u\in V$.
\end{itemize}\par
The next proposition gives the well-posedness of equation (\ref{spde-b}) whose proof can be found in \cite{RSZ}.
\begin{proposition}\label{WP1}
  	Assume (\hyperlink{H1*}{H1*})-(\hyperlink{H5*}{H5*}) hold. Then there is a constant $\varepsilon_0\in(0,\delta]$ such that for any  $0<\varepsilon\leq\varepsilon_0$, equation (\ref{spde-b}) has a unique probabilistic strong solution $u^\varepsilon$ in $L^\alpha([0,T],V)\cap C([0,T],H)$.
\end{proposition}

To state the main result in this section, we consider the skeleton equation as introduced in (\ref{skeleton}): for any $h\in L^2([0,T],U)$,
\begin{numcases}{}\label{skeleton-1}
	dY^h_t=A(t,Y^h_t)dt+B(t,Y^h_t)h(t)dt,
	\text{  $t\in[0,T]$,}\\
	Y^h_0=x\in H\nonumber.
\end{numcases}
The well-posedness of equation (\ref{skeleton-1}) is given in the following proposition, whose proof is deferred to Subsection 4.2.
\begin{proposition}\label{WP2}
    Assume (\hyperlink{\text{H1*}}{\text{H1*}})-(\hyperlink{\text{H5*}}{\text{H5*}}) hold. Then for any initial value $x\in H$, there exists a unique solution $Y^h\in C([0, T], H) \cap L^\alpha([0, T], V )$ to equation (\ref{skeleton-1})  in the sense of Definition \ref{sol} .
\end{proposition}
\par

Now we can state the main result of this section.
\begin{theorem}\label{ldp1}
	Assume (\hyperlink{\text{H1*}}{\text{H1*}})-(\hyperlink{\text{H5*}}{\text{H5*}}) hold. The solutions $\{u^\varepsilon\}_{0<\varepsilon\leq\varepsilon_0}$ to equation (\ref{spde})  satisfy the large deviation principle on $C([0,T],H)$, with the rate function
		\begin{align*}
	I(f)=\inf\limits_{\big\{{h}\in L^2([0,T],U):f=Y^h\big\}}\Big\{\frac{1}{2}\int_{0}^{T}\|{h}(s)\|_U^2ds\Big\},
\end{align*}
with the convention $\inf\{\emptyset\}=\infty$, where $Y^h$ is the unique solution to (\ref{skeleton-1}).
\end{theorem}
\noindent {\bf Proof}.
To prove Theorem \ref{ldp1}, we will apply the weak convergence method outlined in Theorem \ref{criterion}.\par
By a similar argument as in the proof of Theorem \ref{ldp}, there exists a family of maps $\{\mathcal{G}^\varepsilon(\cdot)\}_{0<\varepsilon\leq\varepsilon_0}$ and $\mathcal{G}^0:C([0,T],U_1)\rightarrow C([0,T],H)$ such that for any $U$-cylindrical Brownian motion $W$  and any $h\in S_N$ for some $N<\infty$, $\mathcal{G}^\varepsilon(W)$ and $\mathcal{G}^0(\int_{0}^{\cdot}h(s)ds)$ is the unique solution to equation (\ref{spde-b}) and (\ref{skeleton-1}), respectively.
According to Theorem \ref{criterion}, the rest of the proof will be divided into two parts.
\par
Part one is to verify condition (a) in Theorem \ref{criterion}. This will be done in Subsection 4.3.
Part two is to prove that condition (b) in Theorem \ref{criterion} holds. This is done in Subsection 4.4.
 $\hfill\blacksquare$

\subsection{Skeleton equations}
\ \ \ In this subsection, we prove the well-posedeness of (\ref{skeleton-1}). The proof of uniqueness of the solution to (\ref{skeleton-1}) is the same as that in Subsection 3.2. \par
To demonstrate the existence of solutions, we start by considering the following parabolic equation on the finite-dimensional space $H_n$ $\rm\big($recall $\{e_i\}_{i\geq1}$, $H_n$ and  $P_n$ defined in Section 3 $\rm\big)$:
\begin{numcases}{}
	dY^n_t=P_nA(t,Y^n_t)dt+P_nB(t,Y^n_t)h(t)dt,\ \text{$t\in[0,T]$,}\nonumber\\
	Y^n_0=P_nx\label{approx1}\in H_n.
\end{numcases}
A similar
argument as that in Subsection 3.2 ensures that there exists a unique global solution $Y^n=\big\{Y^n(t),t\in[0,T]\big\}$ to equation $(\ref{approx1})$. Moreover, there exists a constant $C>0$ such that the following inequality holds,
\begin{align}
	\sup_{n\in\mathbb{N}}\big[\sup_{t\in[0,T]}\|Y^n_t\|_H^2+\int_{0}^{T}\|Y^n_s\|_V^\alpha ds\big]\leq C.\label{09192124}
\end{align}
Additionally, using a similar proof as that of Lemma \ref{compact1}, we can establish the following relative compactness.
\begin{lemma}\label{compact3}
	The family $\{Y^n\}_{n\geq1}$ of solutions of equation (\ref{approx1}) is precompact in $C([0,T],V^*)\cap L^2([0,T],H)$.
\end{lemma}
As a corollary of the above lemma, there exists a subsequence (still labeled by $\{Y^n\}_{n\geq1}$) , an element $Y\in C([0,T],V^*)\cap L^2([0,T],H)$, ${\tilde{\mathcal{A}}}\in L^\frac{\alpha}{\alpha-1}([0,T],V^*)$ and an element ${\tilde{\mathcal{B}}}\in L^2([0,T],L_2(U,H))$ such that as $n\rightarrow\infty$,
\begin{align}
	&Y^n\rightarrow Y\text{ weakly in $L^\alpha([0,T],V)$},\\&Y^n\rightarrow Y\text{ in $L^2([0,T],H)\cap C([0,T],V^*)$}\label{02120218}, \\&Y^n\rightarrow Y\text{ in the weak $*$ topology of $L^\infty([0,T],H)$},\\&A(\cdot,Y^n_\cdot)\rightarrow \tilde{\mathcal{A}}\text{ weakly in $L^\frac{\alpha}{\alpha-1}([0,T],V^*)$},\\&P_nB(\cdot,Y^n_\cdot)\rightarrow \tilde{\mathcal{B}}\text{ weakly in $L^2([0,T],L_2(U,H))$}.
\end{align}
For $t\in [0,T]$, let's define
\begin{align*}
	\tilde{Y}_t=x+\int_{0}^{t}\tilde{\mathcal{A}}(s)ds+\int_{0}^{t}\tilde{\mathcal{B}}(s)h(s)ds.
\end{align*}
Using a similar argument as in Section 3, we can show that $\tilde{Y}=Y$, which means that for $a.e.\ t\in [0,T]$,
\begin{align}\label{eq Yt}
	Y_t=x+\int_{0}^{t}\tilde{\mathcal{A}}(s)ds+\int_{0}^{t}\tilde{\mathcal{B}}(s)h(s)ds.
\end{align}
Furthermore, $Y$ is a continuous function valued in $H$. The following lemma shows that $Y$ is a solution to equation (\ref{skeleton-1}).
\begin{lemma}\label{02112314}
	For $\tilde{\mathcal{A}}(\cdot)$, $\tilde{\mathcal{B}}(\cdot)$ and $Y_\cdot$ given above, we have $\tilde{\mathcal{A}}(\cdot)+\tilde{\mathcal{B}}(\cdot)h(\cdot)=A(\cdot,Y_\cdot)+B(\cdot,Y_\cdot)h(\cdot)$, $dt$-$a.s$.
\end{lemma}	
\noindent {\bf Proof}.
The fact that $Y^n_\cdot\rightarrow Y_\cdot$ in $L^2([0,T],H)$ as $n\rightarrow\infty$ implies that, for any $\psi\in L^\infty([0,T],\mathbb{R}_+)$ and  $u_\cdot\in C([0,T],H)\cap L^\alpha([0,T],V)$,
\begin{align}\label{02120137}
	&\int_{0}^{T}\psi(t)\big[\|Y_t\|_H^2-\|x\|_H^2\big]dt\\=&\ \liminf\limits_{n\rightarrow\infty}\int_{0}^{T}\psi(t)\big[\|Y^n_t\|_H^2-\|P_nx\|_H^2\big]dt\nonumber\\=&\ \liminf_{n\rightarrow\infty}\int_{0}^{T}\psi(t)\Bigg\{\int_{0}^{t}2\big\langle A(s,Y^n_s)-A(s,u_s),Y^n_s-u_s\big\rangle ds\nonumber\\&\  +2\int_{0}^{t}\big(B(s,Y^n_s)h(s)-B(s,u_s)h(s),Y^n_s-u_s\big)ds\nonumber\\&\ +2\int_{0}^{t}\langle A(s,Y^n_s),u_s \rangle ds+2\int_{0}^{t}\langle A(s,u_s), Y^n_s\rangle ds-2\int_{0}^{t}\langle A(s,u_s), u_s\rangle ds\nonumber\\&\ +2\int_{0}^{t}\big( B(s,u_s)h(s), Y^n_s-u_s\big) ds+2\int_{0}^{t}\big( B(s,Y^n_s)h(s),u_s\big) ds\Bigg\}dt\nonumber\\ \leq&\  \liminf_{n\rightarrow\infty} \int_{0}^{T}\psi(t)\int_{0}^{t}\big[f(s)+\rho(Y^n_s)+\eta(u_s)+\|h(s)\|_U^2\big]\|Y^n_s-u_s\|_H^2dsdt\nonumber\\&+\ \int_{0}^{T}\psi(t)\Bigg\{\int_{0}^{t}2\big\langle \tilde{\mathcal{A}}(s),u_s\big\rangle ds+2\int_{0}^{t}\big\langle A(s,u_s),Y_s\big\rangle ds-2\int_{0}^{t}\big\langle A(s,u_s),u_s\big\rangle ds\nonumber\\&+\ 2\int_{0}^{t}\big(B(s,u_s)h(s),Y_s-u_s\big)ds+2\int_{0}^{t}\big(\tilde{\mathcal{B}}(s)h(s),u_s\big)ds\Bigg\}dt\label{08121517}.
\end{align}
On the other hand, it follows from (\ref{eq Yt}) that,
\begin{align}\label{02120136}
	&\int_{0}^{T}\psi(t)\big[\|Y_t\|_H^2-\|x\|_H^2\big]dt\\=&\   \int_{0}^{T}\psi(t)\Big\{\int_{0}^{t}\big[2\big\langle \tilde{\mathcal{A}}(s),Y_s\big\rangle+2\big(\tilde{\mathcal{B}}(s)h(s),Y_s\big)\big]ds\Big\}dt\nonumber.
\end{align}
Combining (\ref{08121517}) with (\ref{02120136}) yields
\begin{align}\label{02120204}
	&\int_{0}^{T}\psi(t)\Big\{\int_{0}^{t}\big[2\big\langle \tilde{\mathcal{A}}(s)+\tilde{\mathcal{B}}(s)h(s)-A(s,u_s)-B(s,u_s)h(s),Y_s-u_s\big\rangle\big]ds\Big\}dt\nonumber\\\leq&\ C\liminf_{n\rightarrow\infty}\int_{0}^{T}\psi(t)\int_{0}^{t}\big[f(s)+\rho(Y^n_s)+\eta(u_s)+\|h(s)\|_U^2\big]\|Y^n_s-u_s\|_H^2dsdt.
\end{align}
For $\varepsilon>0$,  $\phi\in L^\infty([0,T],\mathbb{R}_+)$ and $e\in V$ chosen arbitrarily, we substitute $u=Y-\varepsilon\phi e$ in (\ref{02120204}). Next, we divide both sides  of (\ref{02120204}) by $\varepsilon$. Letting $\varepsilon\rightarrow0$ and using (\hyperlink{H1*}{H1*}),  (\hyperlink{H2*}{H2*}), (\ref{09192124}) and the dominated convergence theorem, we obtain the following result
\begin{align}\label{02120341}
	&\int_{0}^{T}\psi(t)\int_{0}^{t}\big\langle \tilde{\mathcal{A}}(s)+\tilde{\mathcal{B}}(s)h(s)-A(s,Y_s)-B(s,Y_s)h(s), e \big\rangle\phi(s)dsdt\leq 0.
\end{align}

The arbitrariness of $e,\phi$ and $\psi$ implies that $\tilde{\mathcal{A}}(\cdot)+\tilde{\mathcal{B}}(\cdot)h(\cdot)=A(\cdot,Y_\cdot)+B(\cdot,Y_\cdot)h(\cdot)$. The proof of Lemma \ref{02112314} is then complete. $\hfill\blacksquare$\par
\begin{remark}\label{08121702}
	By employing similar arguments as in the proof of (\ref{estimate0}), we can demonstrate that for any $N>0$, there exists a constant $C_N>0$ such that
	\begin{align}
	\sup_{h\in S_N}\Big\{\sup_{t\in[0,T]}\|Y^h_t\|_H^2+\int_{0}^{T}\|Y^h_s\|_V^\alpha ds\Big\}\leq C_N<\infty.\label{05210202}
	\end{align}
  Here $Y^h$ is the unique solution to equation (\ref{skeleton-1}).
\end{remark}

\subsection{The proof of the main result: part one}
\ \ \ \ In this subsection, we will verify condition (a) of Theorem \ref{criterion}, which corresponds to part one of the proof of Theorem \ref{ldp1}. Recall the definition of $\mathcal{G}^\varepsilon(\cdot)$ given in the proof of Theorem \ref{ldp1}. By the Girsanov's transformation, for  any $ h^\varepsilon\in \mathcal{A}_N$ and $\varepsilon\in(0,\varepsilon_0]$, the process $X^\varepsilon\dot{=}\ \mathcal{G}^\varepsilon\big(W_\cdot+\frac{1}{\varepsilon}\int_{0}^{\cdot}h^\varepsilon(s)ds\big)$ is  the unique solution of the following stochastic equation,
\begin{numcases}{}
	dX^\varepsilon_t=A(t,X^\varepsilon_t)dt+\varepsilon B(t,X^\varepsilon_t)dW_t+B(t,X^\varepsilon_t)h^\varepsilon(t)dt,\ t\in (0,T],\label{spde3}\\
	X^\varepsilon_0=x\in H.\nonumber
\end{numcases}
We first give the following uniform estimate.
\begin{proposition}\label{02140017}

	For any $q\in[2,\infty)$ and $\{h^\varepsilon\}_{0<\varepsilon\leq\varepsilon_0}\subseteq \mathcal{A}_N$, there exist constants $c_q\in(0,\varepsilon_0]$ and $C_{q,N}>0$ such that
	$$\sup_{0<\varepsilon\leq c_q}E\Big[\big(\int_{0}^{T}\|X^\varepsilon_s\|_V^\alpha ds\big)^\frac{q}{2}+\sup_{t\in[0,T]}\|X_t^\varepsilon\|_H^q\Big]<C_{q,N}<\infty. $$

\end{proposition}
The proof of the this proposition is very similar to that of Proposition \ref{estimate1}. The difference is that we use the new growth conditions on the coefficient $B(t,u)$, (\hyperlink{H3*}{H3*}), (\hyperlink{H5*}{H5*}) instead of the assumptions (\hyperlink{H3}{H3}), (\hyperlink{H5}{H5}). To avoid the repeat, we omit the details.
\vskip 0.5cm
Set
\begin{align}\label{08121637}
	\tau^{\varepsilon,M}\dot{=}\inf\Big\{t\geq0:\|X^\varepsilon_t\|_H\geq M\ \text{or}\ \int_{0}^{t}\|X^\varepsilon_s\|_V^\alpha ds\geq M\Big\}\wedge T.
\end{align}
By Chebyshev's inequality, we have the following corollary.
\begin{corollary}\label{07282053} There exists a constant $c_2\in(0,\varepsilon_0)$ such that as $M\rightarrow\infty$, we have
	$$\sup_{0<\varepsilon\leq c_2}P\Big(\tau^{\varepsilon,M}<T\Big)\rightarrow0.$$
\end{corollary}\par
\noindent\textbf{\textit{Proof of Theorem \ref{ldp1}: part one.}}
\vskip 0.5cm
Let the family $\{h^\varepsilon\}_{0<\varepsilon\leq\varepsilon_0}\subseteq\mathcal{A}_N$ be given. Recall that $X^\varepsilon\dot{=}\ \mathcal{G}^\varepsilon\big(W_\cdot+\frac{1}{\varepsilon}\int_{0}^{\cdot}h^\varepsilon(s)ds\big)$ and $Y^\varepsilon\ \dot{=}\ \mathcal{G}^0(\int_{0}^{\cdot}h^\varepsilon(s)ds)$, which satisfy equations (\ref{spde3}) and (\ref{skeleton-1}) for $\varepsilon\in(0,\varepsilon_0)$, respectively. Therefore we can express the difference between these solutions as
\begin{align*}
	X^\varepsilon_t-Y^\varepsilon_t=&\int_{0}^{t}\Big(A(s,X^\varepsilon_s)-A(s,Y^\varepsilon_s)\Big)ds+\int_{0}^{t}\Big(B(s,X_s^\varepsilon)h^\varepsilon(s)-B(s,Y^\varepsilon_s)h^\varepsilon(s)\Big)ds\\&+\varepsilon\int_{0}^{t}B(s,X_s^\varepsilon)dW_s.
\end{align*}
Applying $\rm It\hat{o}$'s formula and using (\hyperlink{H2*}{H2*}), we have
\begin{align}\label{02132104}
	&\ \|X^\varepsilon_t-Y^\varepsilon_t\|_H^2\nonumber\\=&\ 2\int_{0}^{t}\big\langle A(s,X^\varepsilon_s)-A(s,Y^\varepsilon_s), X^\varepsilon_s-Y^\varepsilon_s \big\rangle ds+2\varepsilon\int_{0}^{t}\Big(B(s,X_s^\varepsilon)dW_s,X^\varepsilon_s-Y^\varepsilon_s\Big)\nonumber\\&\ +2\int_{0}^{t}\Big( B(s,X_s^\varepsilon)h^\varepsilon(s)-B(s,Y^\varepsilon_s)h^\varepsilon(s),X^\varepsilon_s-Y^\varepsilon_s\Big)ds\nonumber\\&\ + \varepsilon^2\int_{0}^{t}\|B(s,X^\varepsilon_s)\|_{L_2}^2ds\nonumber\\\leq&\ \int_{0}^{t}2\big\langle A(s,X^\varepsilon_s)-A(s,Y^\varepsilon_s),X^\varepsilon_s-Y^\varepsilon_s\big\rangle ds+\int_{0}^{t}\|B(s,X^\varepsilon_s)-B(s,Y^\varepsilon_s)\|^2_{L_2}ds\nonumber\\&\ +\int_{0}^{t}\|h^\varepsilon(s)\|_U^2\|X^\varepsilon_s-Y^\varepsilon_s\|_H^2ds+2\varepsilon\int_{0}^{t}\Big(B(s,X^\varepsilon_s)dW_s, X^\varepsilon_s-Y^\varepsilon_s\Big)\nonumber\\&\ +\varepsilon^2\int_{0}^{t}\|B(s,X^\varepsilon_s)\|_{L_2}^2ds\nonumber\\\leq&\ \int_{0}^{t}\Big(f(s)+\|h^\varepsilon(s)\|_U^2+\rho(X^\varepsilon_s)+\eta(Y^\varepsilon_s)\Big)\|X^\varepsilon_s-Y^\varepsilon_s\|_H^2 ds+\varepsilon^2\int_{0}^{t}\|B(s,X^\varepsilon_s)\|_{L_2}^2ds\nonumber\\&\ +2\varepsilon\int_{0}^{t}\Big(B(s,X^\varepsilon_s)dW_s, X^\varepsilon_s-Y^\varepsilon_s\Big).
\end{align}
Replacing $t$ by $t\wedge\tau^{\varepsilon,M}$ in (\ref{02132104}) and using Gronwall's inequality, (\hyperlink{H5*}{H5*}) and (\ref{05210202}), we have
\begin{align}\label{eq zhai 20230712}
\sup_{t\leq\tau^{\varepsilon,M}}\|X_t^\varepsilon-Y_t^\varepsilon\|_H^2\leq C_{M,N}N^{\varepsilon,M}_T,
\end{align}
where $$N^{\varepsilon,M}_T\dot{=}\ \varepsilon^2C_M+2\varepsilon\sup_{t\in[0,\tau^{\varepsilon,M}]}\big|\int_{0}^{t}\Big(B(s,X^\varepsilon_s)dW_s, X^\varepsilon_s-Y^\varepsilon_s\Big)\big|.$$
Applying the BDG inequality, (\hyperlink{H5*}{H5*}) and Remark \ref{08121702}, for $0<\varepsilon\leq c_2$, where $c_2$ is the constant in Corollary \ref{07282053}, we get
\begin{align}\label{eq zhai 2023071200}
	\nonumber E\Big[N^{\varepsilon,M}_T\Big]&\leq\varepsilon^2C_M+2\varepsilon E\Big[\big(\int_{0}^{\tau^{M,\varepsilon}}\|B(s,X^\varepsilon_s)^*(X^\varepsilon_s-Y^\varepsilon_s)\|_U^2ds\big)^\frac{1}{2}\Big]\\&\leq \varepsilon^2C_M+2\varepsilon C_{M,N}.
\end{align}
Then we obtain that for $\forall\ \eta>0$,
\begin{align*}
	&\ P\Big(\sup_{t\leq T}\|X^\varepsilon_t-Y^\varepsilon_t\|_H^2\geq\eta\Big)\\\leq&\  P\Big(\tau^{\varepsilon,M}< T\Big)+P\Big(\tau^{\varepsilon,M}= T,\ \sup_{t\leq T}\|X^\varepsilon_t-Y^\varepsilon_t\|_H^2\geq\eta\Big)\\\leq&\  P\Big(\tau^{\varepsilon,M}<T\Big)+\frac{1}{\eta}E\Big[\sup_{t\leq\tau^{\varepsilon,M}}\|X_t^\varepsilon-Y_t^\varepsilon\|_H^2\Big]\\\leq&\ \sup_{\epsilon\leq{c_2}} P\Big(\tau^{\epsilon,M}<T\Big)+ C_{M,N,\eta}\Big(\varepsilon^2C_{M}+2\varepsilon C_{M,N}\Big).
\end{align*}
Here we have used (\ref{eq zhai 20230712}) and (\ref{eq zhai 2023071200}) to derive the last inequality.
By letting $\varepsilon\rightarrow0$, we have
\begin{align*}
	\lim\limits_{\varepsilon\rightarrow0}P\Big(\sup_{t\leq T}\|X^\varepsilon_t-Y^\varepsilon_t\|_H^2\geq \eta\Big)\leq\sup_{\epsilon\leq{c_2}} P\Big(\tau^{\epsilon,M}<T\Big).{}
\end{align*}
 Next, by letting $M\rightarrow\infty$ and applying Corollary \ref{07282053}, we conclude that
 \begin{align*}
 	\lim\limits_{\varepsilon\rightarrow0}P\Big(\sup_{t\leq T}\|X^\varepsilon_t-Y^\varepsilon_t\|_H^2 \geq\eta\Big)=0.
 \end{align*}
Therefore, the part one of the proof of Theorem \ref{ldp1} is complete. $\hfill\blacksquare$
 \subsection{The proof of the main result: part two}
 \ \ \ \ This subsection is devoted to the verification of (b) in Theorem \ref{criterion}. Recall  that we have defined $Y^n\ \dot{=}\ \mathcal{G}^0\big(\int_{0}^{\cdot}{h}^n(s)ds\big)$ and $Y^h\ \dot{=}\ \mathcal{G}^0\big(\int_{0}^{\cdot}{h}(s)ds\big)$ for $\{h^n,h\}_{n\geq1}\subseteq S_N$. We will show that if $ h^n\rightarrow h$ weakly in $L^2([0,T],U)$ as $n\rightarrow\infty$, then
 \begin{align}\label{eq 20230714}
 	\lim_{n\rightarrow\infty}\sup_{s\in[0,T]}\|Y^n_s-Y^h_s\|_H=0.
 \end{align}

Recall
 (\ref{05210202}) implies that
 \begin{align}\label{estimate3}
 	\sup_{n\geq1}\Big\{\sup_{t\in[0,T]}\|Y^n_t\|_H+\int_{0}^{T}\|Y^n_s\|_V^\alpha ds\Big\}<\infty.
 \end{align}
 Using a similar proof as that of Theorem \ref{compact1}, we can obtain the relative compactness of
 $\{Y^n\}_{n\geq1}$ in $L^2([0,T],H)$.
So, there exists a subsequence, still labeled by  $\{Y^n\}_{n\geq1}$, and an element ${Y}\in L^2([0,T],H)$ such that as $n\rightarrow\infty$,
 \begin{align}
 	&Y^{n}\rightarrow Y\text{  in $L^2([0,T],H)$},\label{05242049}\\&Y^{n}\rightarrow Y\text{  weakly in $L^\alpha([0,T],V)$}.\nonumber
 \end{align}
Moreover, there exists a further subsequence, still labeled by  $\{Y^n\}_{n\geq1}$, such that as $n\rightarrow\infty$,
 \begin{align*}
	&Y^{n}_s\rightarrow Y_s\text{ in $H$ $a.e.\ s\in [0,T]$}.
\end{align*}
(\ref{estimate3}) implies that
 \begin{align}\label{05241755}
	Y\in L^\infty([0,T],H)\cap L^\alpha([0,T],V).
\end{align}
With these preparations, we are now in the position to prove (\ref{eq 20230714}).

Notice that for any $t\in[0,T]$,
$$Y^n_t-Y^h_t=\int_{0}^{t}\Big(A(s,Y^n_s)-A(s,Y^h_s)\Big)ds+\int_{0}^{t}\Big(B(s,Y^n_s)h^n(s)-B(s,Y^h_s)h(s)\Big)ds.$$
Applying the chain rule and using (\hyperlink{H2*}{H2*}), we have
\begin{align}\label{02142055}
	&\ \|Y^n_t-Y^h_t\|_H^2\\=&\ 2\int_{0}^{t}\big\langle A(s,Y^n_s)-A(s,Y^h_s),Y^n_s-Y^h_s\big\rangle ds\nonumber\\&\ +2\int_{0}^{t}\Big(B(s,Y^n_s)h^n(s)-B(s,Y^h_s)h^n(s),Y^n_s-Y^h_s\Big)ds\nonumber\\&\ +2\int_{0}^{t}\Big(B(s,Y^h_s)h^n(s)-B(s,Y^h_s)h(s), Y^n_s-Y^h_s\Big)ds\nonumber\\\leq& \int_{0}^{t}\big[f(s)+\rho(Y^n_s)+\eta(Y^h_s)+\|h^n(s)\|_U^2\big]\|Y^n_s-Y^h_s\|_H^2ds\nonumber\\&\ +2\int_{0}^{t}\Big(B(s,Y^h_s)h^n(s)-B(s,Y^h_s)h(s),Y^n_s-Y^h_s\Big)ds\nonumber.
\end{align}
Using Gronwall's inequality and (\ref{estimate3}), we arrive at
\begin{align}\label{05242109}
	\sup_{t\in[0,T]}\|Y^n_t-Y^h_t\|_H^2\leq C_N\sup_{t\in[0,T]}|Z_n(t)|.
\end{align}
Here
\begin{align*}
Z_{n}(t)\ \dot{=}\int_{0}^{t}\Big(B(s,Y^h_s)\big(h^n(s)-h(s)\big), Y^n_s-Y^h_s\Big)ds.
\end{align*}

To verify (\ref{eq 20230714}), we need to prove that $Z_n\rightarrow0$ in $C([0,T],\mathbb{R})$ as $n\rightarrow\infty$.
By (\ref{estimate3}), for any $t,s\in[0,T]$ with $s<t$, we have
\begin{align*}
|Z_n(t)-Z_n(s)|&\leq\int_{s}^{t}\big|\Big( B(r,Y^h_r)\big(h^n(r)-h(r)\big),Y^n_r-Y^h_r\Big)\big|dr\\&\leq \  C\Big(\int_{s}^{t}\|B(r,Y^h_r)\|_{L_2}^2dr\Big)^\frac{1}{2}N^\frac{1}{2}.
\end{align*}
Thus, $\{Z_n\}_{n\geq1}$ is an equi-continuous family in $C([0,T],\mathbb{R})$. Next, we will show that for any $t\in[0,T]$, $Z_n(t)\rightarrow 0$. For any $t\in[0,T]$,
\begin{align*}
	Z_{n}(t)\ &{=}\int_{0}^{t}\Big(B(s,Y^h_s)\big(h^n(s)-h(s)\big), Y^n_s-Y_s\Big)ds\\&\ \ \  +\int_{0}^{t}\Big(B(s,Y^h_s)\big(h^n(s)-h(s)\big), Y_s-Y^h_s\Big)ds\\&=\ \text{I}^n(t) + \text{II}^n(t).
\end{align*}
Since
$ h^n\rightarrow h$ weakly in $L^2([0,T],U)$ as $n\rightarrow\infty$
and $B(\cdot,Y^h_\cdot)^*(Y_\cdot-Y^h_\cdot)\in L^2([0,T],U)$, then for any $t\in[0,T]$, $$\lim_{n\rightarrow\infty}\text{II}^n(t)= 0.$$
 For any $M>0$, set
\begin{align*}
	\Gamma_M\dot{=}\big\{t\in[0,T]:\|B(t,Y^h_t)\|_{L_2}>M\big\}.
\end{align*}
Using Chebyshev's inequality and (\ref{estimate3}), we can estimate the Lebesgue measure $|\Gamma_M|$ of the set $\Gamma_M$ as follows
\begin{align}\label{05242041}
|\Gamma_M|\leq \frac{C_N}{M^2}.
\end{align}
By (\ref{estimate3}) again, we have
\begin{align*}
\text{I}^n(t)&\leq C\int_{\Gamma_M}\|h^n(s)-h(s)\|_H\|B(s,Y^h_s)\|_{L_2}ds+M\int_{0}^{T}\|h^n(s)-h(s)\|_U\|Y_s-Y^n_s\|_Hds\\&\leq C_N\Big(\int_{\Gamma_M}\|B(s,Y^h_s)\|_{L_2}^2ds\Big)^\frac{1}{2}+MC_N\Big(\int_{0}^{T}\|Y_s-Y^n_s\|_H^2ds\Big)^\frac{1}{2}.
\end{align*}
Letting $n\rightarrow+\infty$ and utilizing $(\ref{05242049})$, we obtain
$$\lim\limits_{n\rightarrow\infty}\text{I}^n(t)\leq C_N\Big(\int_{\Gamma_M}\|B(s,Y^h_s)\|_{L_2}^2ds\Big)^\frac{1}{2}.$$
Letting $M\rightarrow\infty$, and applying (\hyperlink{H5*}{H5*}) and $(\ref{05242041})$, we can assert that
$$\lim\limits_{n\rightarrow\infty}\text{I}^n(t)=0.$$
 Therefore, for any $t\in[0,T]$, $Z_n(t)\rightarrow 0$. Combining this with the equi-continuity of $\{Z_n\}_{n\geq0}$, the Arzela-Ascoli's theorem shows that $Z_n\rightarrow0$ in $C([0,T],\mathbb{R})$ as $n\rightarrow\infty$. Due to the uniqueness of the solutions of the skeleton equation, the limit (\ref{eq 20230714}) holds for the whole sequence. Thus  the part two of the proof of Theorem \ref{ldp1} is completed.
$\hfill\blacksquare$


\section{Applications}\label{sec 5}
\setcounter{equation}{0}
\ \ \ \ In this section, we will provide examples of different types of SPDEs that satisfy our framework in Parts A or B. \par In fact, our framework can cover a great number of interesting examples.  Specificially, all the examples considered in \cite{LR2,L,RSZ} fulfill our framework. For instance, the 2D Navier-Stokes equations, the 3D tamed Navier-Stokes equations, p-Laplacian equations, 1D Burgers equations, porous media equations, Allen-Cahn equations, Cahn-Hilliard equations, fast-diffusion equations, 3D Leray-$\alpha$ model, 2D Boussinesq system, 2D MHD equations, 2D Boussinesq model for the B\'{e}nard equations, some shell models of turbulence (GOY, Sabra, dyadic), power law fluids, the Ladyzhenskaya model and the Kuramoto-Sivashinsky equations are included.\par Moreover, our main results can also be applied to stochastic quasilinear PDEs and $p$-Laplace equations with nonlinear transport type noise, whose LDP  have not been covered by any other framework or proven in existing literature.
 Our results can also improve some existing results on the large deviation principle of SPDEs, including the convection diffusion equation, liquid crystal system, Cahn-Hilliard-Navier-Stokes equation, Allen-Cahn-Navier-Stokes equation and in particlar, the strong solution for the 2D Allen-Cahn equation and the two-dimensional Navier-Stokes equation.\par
    The following examples will illustrate that the conditions (\hyperlink{H1}{H1})-(\hyperlink{H5}{H5}) or (\hyperlink{H1*}{H1*})-(\hyperlink{H5*}{H5*}) are quite general and can be verified easily. The first three examples, namely, the quasilinear equations, the convection diffusion equation, and the 2D Liquid crystal equation, are applications of the results in Part A. The last two examples, the $p$-Laplacian equation and the 2D Navier-Stokes equation perturbed by gradient-dependent noise, are applications of the results in Part B.
\begin{example}\label{qsl}(Quasilinear SPDEs). Let $\mathcal{D}$ be a bounded domain in $\mathbb{R}^d$ with smooth boundary $\partial{\mathcal{D}}$. We consider the following stochastic  quasilinear partial differential equation:
	\begin{numcases}{}
	\nonumber	\partial_tu^\varepsilon(t,x)=\nabla\cdot a\big(t,x,u^\varepsilon(t,x),\nabla u^\varepsilon(t,x)\big)-a_0\big(t,x,u^\varepsilon(t,x),\nabla u^\varepsilon(t,x)\big),\\\qquad\qquad\ \ \ \nonumber+\ \varepsilon B\big(t,u^\varepsilon(t,x)\big)\dot{W}_t,\text{  \ \  for $(t,x)\in$ $(0,T]\times\mathcal{D}$,}\\
	\nonumber	u^\varepsilon(t,x)=0,\qquad\qquad \text{ for $(t,x)\in(0,T]\times\partial\mathcal{D}$},\\
	\nonumber	u^\varepsilon(0,x)=u_0(x),\qquad\qquad x\in\mathcal{D}.
	\end{numcases}
	Here $u^\varepsilon:[0,T]\times\mathcal{D}\rightarrow\mathbb{R}$ represents the solution, the vector $\nabla u^\varepsilon(t,x)=\big(\partial_iu^\varepsilon(t,x)\big)_{i=1}^d$ is the gradient of $u^\varepsilon$ with respect to the spatial variable $x$, $a=(a_1,a_2,\cdots,a_d)$ is a vector with $a_i:[0,T]\times\mathcal{D}\times\mathbb{R}\times\mathbb{R}^d\rightarrow\mathbb{R}$ for each $i=1,\cdots,n$.\par
	We assume that $a_i,i=0,1,\cdots,d$, satisfy the following conditions:
	there exists a constant $\alpha>1$ for $d=1,2$ and $\alpha\geq\frac{2d}{d+2}$ for $d\geq3$, such that\par
	\begin{itemize}
		\item [\hypertarget{S1}{{\bf (S1)}}] For each $i=1,\cdots,d$, $a_i$ satisfies the Carath\'{e}odory condition: for $a.e.$ fixed $(t,x)\in [0,T]\times\mathcal{D}$, $a_i(t,x,u,z)$ is continuous in $(u,z)\in\mathbb{R}\times{\mathbb{R}^d}$, and for each fixed $(u,z)\in\mathbb{R}\times\mathbb{R}^d$, $a_i(t,x,u,z)$ is measurable with respect to $(t,x)\in[0,T]\times\mathcal{D}$.
		\item [\hypertarget{S2}{{\bf (S2)}}] There exist nonnegative constants $c_1$ and $c_2$ and a function $f_1\in L^\frac{\alpha}{\alpha-1}([0,T]\times\mathcal{D},\mathbb{R}_+)$ such that for $a.e.$ $(t,x)\in[0,T]\times\mathcal{D}$ and all $(u,z)\in\mathbb{R}\times\mathbb{R}^d,i=1,\cdots,d,$
		$|a_i(t,x,u,z)|\leq c_1|z|^{\alpha-1}+c_2|u|^\frac{(\alpha-1)(d+2)}{d}+f_1(t,x).$
		\item [\hypertarget{S3}{{\bf (S3)}}] There exists constant $c_3>0$, $c_4\geq 0$, and a function $f_2\in L^1([0,T]\times\mathcal{D},\mathbb{R}_+)$ such that for $a.e.$ $(t,x)\in[0,T]\times\mathcal{D}$ and all $(u,z)\in\mathbb{R}\times\mathbb{R}^d$,
		$$\sum_{i=1}^{d}a_i(t,x,u,z)z_i+a_0(t,x,u,z)u\geq c_3|z|^\alpha-c_4|u|^2-f_2(t,x).$$
		\item [\hypertarget{S4}{{\bf (S4)}}] For $a.e.$ $(t,x)\in[0,T]\times\mathcal{D}$, all $u\in\mathbb{R}$ and $z,\tilde{z}\in\mathbb{R}^d$ such that $z\neq\tilde{z}$,
		$$\sum_{i=1}^{d}\big[a_i(t,x,u,z)-a_i(t,x,u,\tilde{z})\big](z_i-\tilde{z}_i)>0.$$
		And for $a.e.$ $(t,x)\in[0,T]\times\mathcal{D}$, and any $M>0$,
		$$\lim_{|z|\rightarrow\infty}\frac{\sup_{|u|\leq M}\sum_{i=1}^{d}a_i(t,x,u,z)z_i}{|z|+|z|^{\alpha-1}}=\infty. $$
		\item [\hypertarget{S5}{{\bf (S5)}}] Let $0\leq\gamma\leq\alpha(1+\frac{2}{d})-2$ and $f_3\in L^1([0,T],\mathbb{R}_+).	$ There exists a constant $c>0$ such that for $a.e.$ $(t,x)\in[0,T]\times\mathcal{D}$ and all  $u,\tilde{u}\in\mathbb{R}$ and $z,\tilde{z}\in\mathbb{R}^d$,
		\begin{align*}
			&\sum_{i=1}^{d}\big[a_i(t,x,u,z)-a_i(t,x,\tilde{u},\tilde{z})\big](z_i-\tilde{z}_i)\\&+\big[a_0(t,x,u,z)-a_0(t,x,\tilde{u},\tilde{z})\big](u-\tilde{u})\geq -c\big(f_3(t)+|u|^\gamma+|\tilde{u}|^\gamma\big)|u-\tilde{u}|^2.\nonumber
		\end{align*}
	\end{itemize}
	In this case, $H\dot{=}\ L^2(\mathcal{D})$ and $V\dot{=}\ W^{1,\alpha}_0(\mathcal{D})$, the Sobolev space with zero trace. For $u,v\in V$, the operator $A$ is defined as follows
	\begin{align*}
		\langle A(t,u),v\rangle=-\int_{0}\Big\{\sum_{i=1}^{d}a_i(t,x,u(x),\nabla u(x))\partial_i v(x)+a_0(t,x,u(x),\nabla u(x))v(x)\Big\}dx.
	\end{align*}
	For the verification of (\hyperlink{H1}{H1})-(\hyperlink{H5}{H5}) for stochastic quasilinear equations, we refer the readers to Section 4 of \cite{RSZ} for a detailed explanation. Therefore, under appropriate conditions on the diffusion term $B(t,\cdot)$ (e.g, $B(t,\cdot)$ satisfies (\hyperlink{H2}{H2}), (\hyperlink{H3}{H3}) and (\hyperlink{H5}{H5}) with $A(t,\cdot)$ provided above), the LDP holds. In particular, if $B(t,\cdot)$ satisfies the linear growth and global Lipschitz condition, then the LDP holds.
\end{example}
\begin{remark}
The large deviation principle for quasilinear stochastic partial differential equations was studied in \cite{MSZ}. However,  the assumptions here on the coefficients are much weaker and more general.
\end{remark}
\begin{example}(Convection diffusion equation).
	The convection diffusion equation describes physical phenomena where particles, energy, or other physical quantities are transferred inside a physical system through two processes: diffusion and convection. This equation has significant applications in fluid dynamics, heat transfer, and mass transfer. The stochastic  convection diffusion equation is given by
	\begin{numcases}{}\label{cde}
		\partial_t u^\varepsilon=\nabla\cdot\big[a(u^\varepsilon)\nabla u^\varepsilon+b(u^\varepsilon)\big]+\varepsilon\sigma(u^\varepsilon)\dot{W}_t,\ on\ (0,T]\times \mathbb{T}^d,\\
		u^\varepsilon(0)=u_0,\nonumber
	\end{numcases}
	where $\mathbb{T}^d\ \dot{=}\ (\mathbb{R}/\mathbb{Z})^d$, the function $u^\varepsilon:[0,T]\times\mathbb{T}^d\rightarrow\mathbb{R}$ represents the solution, the vector $b=\ (b_1,\cdots,b_d):\mathbb{R}\rightarrow\mathbb{R}^d$ denotes the flux function, the matrix $a=(a_{ij}):\mathbb{R}^d\rightarrow\mathcal{M}_{d\times d}$ is called the diffusion matrix. Here $\mathcal{M}_{d\times d}$ is the set of all $d\times d$-dim matrices.  We assume that the flux function  b is globally Lipschitz continuous and  the diffusion matrix a is bounded, global Lispchitz continuous. Moreover, the diffusion matrix $a$ is uniformly elliptic, that is, there exist constants $\delta>0$ and $M>0$ such that for any $u\in\mathbb{R}$ and $z\in\mathbb{R}^d$,
	$$\delta|z|^2\leq\big\langle a(u)z,z \big\rangle\leq M|z|^2.$$
	In this case, we define $H\dot{=}\ L^2(\mathbb{T}^d)$ and $V\dot{=}\ W^{1,2}(\mathbb{T}^d)$. The operator $A$ is given by
	$$\big\langle A(u),v\big\rangle \dot{=} -\int_{\mathbb{T}^d}\big\langle a(u(x))\nabla u(x)+b(u(x)),\nabla v(x)\big\rangle dx \text{ for any $u,v\in V$.}$$
	Moreover, for each $u\in H$, $\sigma(u):U\rightarrow H$ is a map defined as $$\big[\sigma(u)\bar{e}_k\big](\cdot)\dot{=}\ \sigma_k\big(u(\cdot)\big).$$ Here $U$ is a separable Hilbert space with orthonormal basis $(\bar{e}_k)_{k\geq1} $, $W$ is an $U$-cylindrical Wiener process and $\sigma_k(\cdot):\mathbb{R}\rightarrow\mathbb{R}$ are real-valued functions. Suppose that $\sigma$ satisfies the global Lipschitz condition and the linear growth condition
	$$\sum_{i=1}^{\infty}|\sigma_k(y_1)-\sigma_k(y_2)|^2\leq C|y_1-y_2|^2,\text{\ $\forall\ y_1,y_2\in\mathbb{R}$}.$$
	$$\sum_{i=1}^{\infty}|\sigma_k(y)|^2\leq C(1+|y|^2),\ \forall y\in\mathbb{R}.$$
	For the verification of  (\hyperlink{H1}{H1}), (\hyperlink{H2'}{H2'}), (\hyperlink{H3}{H3})- (\hyperlink{H5}{H5}) and the well-posedness of the stochastic equation (\ref{cde}), we refer the readers to Example 4.2 in \cite{RSZ} for details. The existence of the skeleton equation can be obtained by Proposition $\ref{WP0}$, while the uniqueness can be proven  using the same technique as in \cite{DZZ}. Although (\hyperlink{H2}{H2}) does not hold, which has been used in the verification of LDP in Section 3.3 and 3.4, we can still show LDP for equation (\ref{cde}), using same techniques as in \cite{DZZ}.
\end{example}
\begin{remark}
	In \cite{DZZ}, the LDP of equation (\ref{cde}) was established under the assumption that the flux function $b=(b_1,\cdots,b_d):\mathbb{R}\rightarrow\mathbb{R}^d$ and the diffusion matrix $a=(a_{ij}):\mathbb{R}^d\rightarrow\mathcal{M}_{d\times d}$ are continuously differentiable and lipschitz continuous. We do not need the assumptions on differentiability in this example. Using the same argument in this example, we can also prove the LDP of the perturbed p-Laplacian equation considered in \cite{VZ}.
\end{remark}

\begin{example}(2D Liquid crystal equation)\label{lce1}
	The liquid crystal equation demonstrates the temporal evolution of the hydrodynamics of liquid crystals. It is a simplified form of the Ericksen-Leslie system with Ginzburg-Landau approximation, which was established in \cite{LL}. We consider the following stochastic version of the 2D liquid crystal equation.
	\begin{numcases}{}\label{lce}
	\nonumber	\partial_t u^\varepsilon=\Delta u^\varepsilon-(u^\varepsilon\cdot\nabla) u^\varepsilon-\nabla p^\varepsilon-\nabla\cdot (\nabla n^\varepsilon\otimes\nabla n^\varepsilon)+\varepsilon\sigma_1(t,u^\varepsilon,n^\varepsilon)\dot{W}^1_t\text{, in $(0,T)\times\mathcal{D}$},\\
		\partial_t n^\varepsilon=\Delta n^\varepsilon-(u^\varepsilon\cdot\nabla)n^\varepsilon-\Phi(n^\varepsilon)+\varepsilon\sigma_2(t,u^\varepsilon,n^\varepsilon)\dot{W}^2_t,\text{ in $(0,T)\times\mathcal{D}$}\nonumber ,\\
		\nabla\cdot u^\varepsilon=0,\nonumber\text{ in $(0,T)\times\mathcal{D}$},\\
		u^\varepsilon=0\text{ and }\frac{\partial n^\varepsilon}{\partial \nu}=0, \text{ on $(0,T)\times\partial{\mathcal{D}}$},\nonumber\\
		u^\varepsilon(0)=\ u_0\text{ and }\ n^\varepsilon(0)=\ n_0,\nonumber \text{ in $\mathcal{D}$}.
	\end{numcases}

Here, $\mathcal{D}$ represents a bounded domain in $\mathbb{R}^2$ with smooth boundary $\partial \mathcal{D}$. The functions $u^\varepsilon:[0,T]\times\mathcal{D}\rightarrow\mathbb{R}^2$ , $p^\varepsilon:[0,T]\times\mathcal{D}\rightarrow\mathbb{R}$, $n^\varepsilon:[0,T]\times\mathcal{D}\rightarrow\mathbb{R}^3$ represent the velocity, pressure, and director field of the liquid crystal molecules, respectively. Additionally, $\nu$ represents the outward unit normal vector on $\partial \mathcal{D}$. $W_1$ and $W_2$ are two independent cylindrical Wiener processes. By the symbol $\nabla n\otimes\nabla n$, we mean a $2\times2$ matrix with entries defined by
$$(\nabla n\otimes\nabla n)_{ij}=\sum_{k=1}^{3}(\partial_i n_k)(\partial_j n_k). $$
We assume that $\Phi:\mathbb{R}^3\rightarrow\mathbb{R}^3$ is of the form
$$\Phi(n)=\varphi(|n|^2)n=\Big(\sum_{i=0}^{k}a_i|n|^{2i}\Big)n,$$
where $\varphi:[0,\infty)\rightarrow\mathbb{R}$ is a $k$-th real-polynomial  and $a_k>0$. Let $V\dot{=}\ \Big\{u\in H^1(\mathcal{D})^2:\nabla\cdot u=0, u|_{\partial\mathcal{D}}=0\Big\}$ and denote by $H$ the closure of $V$ under the $L^2$-norm $\|u\|_H^2\dot{=}\int_{\mathcal{D}}|u(x)|^2dx.$ Furthermore, we define
$$\mathbb{H}\ \dot{ =}\ H\times[H^1(\mathcal{D})]^3,\ \mathbb{V}\ \dot{=}\  V\times\ \Big\{n\in\big[ H^2(\mathcal{D})\big]^3:\frac{\partial n}{\partial \nu}=\ 0\Big\},$$
where the norms in $\mathbb{H}$ and $\mathbb{V}$ are separately denoted by
$$\|X\|_\mathbb{H}^2\ \dot{=}\ \|u\|_H^2+\|n\|_{H^1}^2,\ \|X\|_{\mathbb{V}}^2\ \dot{=}\ \|u\|_V^2+\|n\|_{H^2}^2,$$
for $X\ =\ (u,n)\in \mathbb{V}.$ It can be seen that we have a Gelfand triple $\mathbb{V}\subseteq\mathbb{H}\subseteq\mathbb{V}^*$ and a compact embedding $\mathbb{V}\subseteq\mathbb{H}$.
The operator $A(\cdot):\mathbb{V}\rightarrow\mathbb{V}^*$ is defined as follows
$$
A(X)\dot{=}\left(
\begin{array}{lcl}
	&P_H[\Delta u-(u\cdot\nabla) u-\nabla n\cdot\Delta n ] \\&\Delta n-(u\cdot\nabla)n-\varphi(n)
\end{array} \right),\text{ for $X=(u,n)\in \mathbb{V}$,}
$$
where $P_H$ represents the usual Helmholtz-Leray projection.\par
\vskip 0.25cm
 With suitable assumptions on the diffusion coefficient $\sigma_1$ and $\sigma_2$ (e.g. (\hyperlink{H2}{H2})(\hyperlink{H3}{H3})(\hyperlink{H5}{H5}) with $A(X)$ defined above) , we can verify conditions (\hyperlink{H1}{H1})-(\hyperlink{H5}{H5}), using the estimates given in Example 4.5 in \cite{RSZ}. Therefore, the small noise LDP holds for equation (\ref{lce}). In particular, for $\sigma_1$ and $\sigma_2$ that are global Lipschitz and of linear growth in $H$, the LDP holds.
\end{example}
\begin{remark}\label{08060102}
	The large deviation principle of the 2D stochastic liquid crystal equation has been established in \cite{BMP} and \cite{ZZ}. However, we improve their results by allowing for more general diffusion coefficients $\sigma_1$ and $\sigma_2$. This framework can also extend the results on the Cahn-Hilliard-Navier-Stokes equation(whose LDP was shown in \cite{DM2} and \cite{ZH}) and the Allen-Cahn-Navier-Stokes equation(whose LDP was shown in \cite{M} and \cite{DM}) by allowing for more general diffusion coefficients. We like to mention that the strong solution of 2D Allen-Cahn equation in the sense of PDE is covered by this framework, which can not be addressed by the results in \cite{LR2}.
\end{remark}\par
Both of the following two examples illustrate the large deviation principle for SPDEs with gradient-dependent noise, where the results in part B apply.
\begin{example}(p-Laplace equation)
	We consider the following stochastic  p-Laplace equation with nonlinear transport type noise.
	\begin{numcases}{}\label{nlple}
		\partial_tu^\varepsilon=\nabla\cdot \Big(|\nabla u^\varepsilon|^{p-2}\nabla u^\varepsilon\Big)+\varepsilon|\nabla u^\varepsilon|^\frac{p}{2}\dot{\beta}_t,\text{ in $(0,T]\times\mathcal{D}$},\\
		\nonumber u^\varepsilon(t,x)=0,\qquad\qquad \text{ for $(t,x)\in\ (0,T]\times\partial\mathcal{D}$},\\
		\nonumber u^\varepsilon(0,x)=u_0(x),\qquad\qquad x\in\mathcal{D}.
	\end{numcases}
Here, $\mathcal{D}$ is a bounded domain in $\mathbb{R}^d$ with smooth boundary. $H=L^2(\mathcal{D})$, $V=W_0^{1,p}(\mathcal{D})$ with $p\geq2$, $U=\mathbb{R}$ and $\beta_\cdot$ denotes an 1-dimensional standard Brownian motion. The notation $|\cdot|$ refers to the Euclidean norm in $\mathbb{R}^d$, while $\cdot$ indicates the usual inner product in $\mathbb{R}^d$. For $u,v\in V$, the operator $A(u)$ is determined by the following inner product,
\begin{align*}
\langle A(u),v\rangle =-\int_{\mathcal{D}}	|\nabla u(x)|^{p-2} \nabla u(x) \cdot \nabla v(x) dx
\end{align*}	
As in Example \ref{qsl}, we can verify that conditions (\hyperlink{H1*}{H1*}), (\hyperlink{H4*}{H4*}) are satisfied with $\alpha=p$. Additionally, the following inequality hold
\begin{align*}
	\|B(t,u)\|^2_{L_2}\leq \|\nabla u\|_{L^p(\mathcal{D})}^p,
\end{align*}
This implies that the condition (\hyperlink{H5*}{H5*}) is satisfied.
Furthermore, according to Lemma 2.2 in \cite{NT}, the following inequalities can be obtained for any $u, v\in V$,
\begin{align*}
	&\langle A(u)-A(v), u-v\rangle\leq -2^{2-p}\|u-v\|_V^p \\&
	\|B(t,u)-B(t,v)\|_{L_2}^2=\int_{\mathcal{D}}\big||\nabla u(x)|^\frac{p}{2}-|\nabla v(x)|^\frac{p}{2}|^2dx\leq\int_{\mathcal{D}}\big|\nabla u(x)-\nabla v(x)|^pdx.
\end{align*}
Hence, condition (\hyperlink{H2*}{H2*}) is satisfied. Moreover, by setting $v=0$ in (\hyperlink{H2*}{H2*}), we can also obtain (\hyperlink{H3*}{H3*}). Therefore, for equation (\ref{nlple}), the large deviation principle stated in Theorem \ref{ldp1} holds.$\hfill\blacksquare$
\begin{remark}
	Here we mention that the large deviation principle of the equation considered above has not been studied by any other paper and can not be covered by any other existing framework.
\end{remark}
\end{example}
Next, we consider the 2D Navier-Stokes equation with transport type noise, which generalizes the result in \cite{SS}.
\begin{example}(2D Navier-Stokes equation)\label{nse1}
	The Navier-Stokes equation is an important model for atmosphere and ocean dynamics, water flow, and other viscous flow.
	The stochastic Navier-Stokes equation with transport-type noise has been studied by many people recently. Examples of these studies include  \cite{SS,CM, FL, GL, PYZ}, which analyzes various properties of the stochastic Navier-Stokes equation with transport-type noise. The two-dimensional model is given by
	\begin{numcases}{}\label{nse}
	\partial_t u^\varepsilon=\Delta u^\varepsilon-(u^\varepsilon\cdot\nabla) u^\varepsilon-\nabla p^\varepsilon+\varepsilon\sigma(t,u^\varepsilon,\nabla u^\varepsilon )\dot{W}_t,\text{ in $(0,T)\times\mathcal{D}$},\\
	\nabla\cdot u^\varepsilon=0,\nonumber\text{ in $(0,T)\times\mathcal{D}$},\\
	u^\varepsilon=0, \text{ on $(0,T)\times\partial{\mathcal{D}}$},\nonumber\\
	u^\varepsilon(0)=u_0,\ \nonumber \text{ in $\mathcal{D}$}.
\end{numcases}
Here, $\mathcal{D}$ is a bounded domain of $\mathbb{R}^2$ with smooth boundary $\partial\mathcal{D}$. The velocity, denoted as  $u^\varepsilon:[0,T]\times\mathcal{D}\rightarrow \mathbb{R}^2$, and the pressure term $p^\varepsilon:[0,T]\times\mathcal{D}\rightarrow\mathbb{R}$, are defined on this domain. $W$ is an $U$-cylindrical Brownian motion and $\sigma(t,u, \nabla u)\in L_2(U,[L^2(\mathcal{D})]^2)$.\par Similar to Example \ref{lce}, let's define the following spaces $$
	H\ \dot{=}\ \Big\{u\in \big[L^2(\mathcal{D})\big]^2:u=0 \ on\  \partial\mathcal{D}\  and\ \nabla\cdot u=0\ in\ \mathcal{D} \Big\}\ $$ and$$ \ V\ \dot{=}\ \Big\{u\in \big[H_0^1(\mathcal{D})\big]^2:\ \nabla\cdot u=0\ in\ \mathcal{D} \Big\}.
$$
We let $\mathcal{P}:\big[L^2(\mathcal{D})\big]^2\rightarrow H$ to be the usual Helmholtz-Leray projection. Hence, the diffusion coefficient can be interpreted as $$
	B(t,u)\ \dot{=}\ \mathcal{P}\sigma(t,u,\nabla u)\in L_2(U,H),\text{ for $u\in V$.}$$
We define the operator $A:V\rightarrow V^*$ as
\begin{align*}
A(u)\ \dot{=}\ \mathcal{P}\big[\Delta u-(u\cdot\nabla) u\big]\ ,\text{ for $u\in V$}.
\end{align*}
Using a similar estimate as in Example \ref{lce}, we see that under suitable conditions on $B(t,\cdot)$(e.g. $B(t,\cdot)$ satisfies (\hyperlink{H2*}{H2*}), (\hyperlink{H3*}{H3*}) and (\hyperlink{H5*}{H5*}) with $A(\cdot)$ constructed above),  (\hyperlink{H1*}{H1*})-(\hyperlink{H5*}{H5*}) hold, which implies the LDP. \par In particular, the LDP holds for the stochastic 2D Navier-Stokes equation when $\sigma(t,\cdot,\cdot)$ satisfies the following Lipschitz and linear growth conditions in the $V$-norm. Namely, there exists $f\in L^1([0,T],\mathbb{R}_+)$ such that for $u,v\in V$,
\begin{align}\label{lip}
	&\|\sigma(t,u,\nabla u)-\sigma(t,v, \nabla v)\|^2_{L_2}\leq C\|u-v\|^2_V+f(t)\|u-v\|^2_H,\\&
	\|\sigma(t,u, \nabla u)\|_{L_2}^2\leq f(t)(1+\|u\|_H^2)+C\|u\|_V^2\label{lg}.
\end{align}	
\vskip -1.5cm
\begin{remark}
This result extends the results in \cite{SS} by allowing more general diffusion coefficient. We like to mention that a similar extension can be applied to the 2D liquid crystal system considered in Example \ref{lce1}.
\end{remark}
\end{example}
\vskip 0.4cm
{\bf Acknowledgement.} 
This work is partially supported by the National Key R\&D Program of China (No. 2022YFA1006001), the National Natural Science Foundation of China (No. 12131019, No. 11721101, No. 12001516), the Fundamental Research Funds for the Central Universities (No. WK3470000031, No. WK3470000024, No. WK3470000016), and the School Start-up Fund (USTC) KY0010000036.

\end{document}